\documentclass[review]{elsarticle}
\usepackage{hyperref}
\usepackage{amsmath, amssymb}
\usepackage{mathtools}
\usepackage{xcolor}
\usepackage{bm}
\usepackage{mathdots}
\usepackage{geometry}
\usepackage{amsthm}
\usepackage{graphicx}
\usepackage{booktabs}
\usepackage{caption}
\usepackage{subcaption}
\usepackage{float}

\theoremstyle{plain}
\newtheorem{theorem}{Theorem}[section]
\newtheorem{lemma}[theorem]{Lemma}

\newtheorem{Remark}[theorem]{Remark}
\newtheorem{Definition}[theorem]{Definition}
\journal{Journal of \LaTeX\ Templates}

\begin{document}

\begin{frontmatter}

\title{Block preconditioning for all-at-once variable-coefficient fractional evolution equations via the GLT analysis}
\author[label1]{Muhammad Faisal Khan$^*$}\ead{mfkhan@uninsubria.it}
\author[label1,label4]{Stefano Serra-Capizzano}\ead{s.serracapizzano@uninsubria.it}
\address[label1]{Department of Science and High Technology, University of Insubria, Como, Italy}
\address[label4]{Department of Information Technology, Uppsala University, Uppsala, Sweden}
\begin{abstract}
We study a class of nonlocal evolutionary partial differential equations with weakly singular temporal kernel and spatially variable diffusion coefficient. The model is posed on $\Omega \subset \mathbb{R}$, and involves a left-sided Riemann--Liouville fractional derivative in space multiplied by a variable coefficient $a(x)$. The temporal derivative is approximated by an $L1$ type scheme, while the spatial operator is discretized by finite difference techniques, resulting in large scale all at once linear systems with a twolevel Toeplitz like structure. We develop and analyze a block lower triangular strategy that mimics the structure of the coefficient matrix while simplifying its components for computational efficiency. The analysis is carried out at the level of matrix sequences by means of generalized locally Toeplitz (GLT) theory. Within this framework, we characterize the asymptotic spectral distribution of the discretized operators and use the associated GLT symbol to guide the construction of the structured approximation. Numerical experiments using the GMRES solver demonstrate that the proposed preconditioning strategy significantly improves convergence rates, robustness, and scalability for large-scale problems. Open problems and possible extensions are briefly discussed at the end of the present work.
\end{abstract}


\end{frontmatter}

\section{Introduction}
Consider the following nonlocal evolutionary fractional differential equation (FDE) with a weakly singular temporal kernel
\begin{equation}\label{1d}
\left\{
\begin{array}{l}
\displaystyle
\frac{1}{\Gamma(1-\beta)}
\int_{0}^{t}
\frac{\partial u(\bm{x},s)}{\partial s}\,(t-s)^{-\beta}\,ds
=
\mathcal{L}u(\bm{x},t) + f(\bm{x},t),
\quad \bm{x}\in \Omega \subset \mathbb{R}^d,\ t\in(0,T], \\[2mm]
\displaystyle
u(\bm{x},t)=0,
\quad \bm{x}\in \partial \Omega,\ t\in(0,T], \\[2mm]
\displaystyle
u(\bm{x},0)=\psi(\bm{x}),
\quad \bm{x}\in \Omega,
\end{array}
\right.
\end{equation}
where $\Gamma(\cdot)$ denotes the Gamma function, $0<\beta<1$, the boundary of $\Omega$ is $\partial \Omega$; $\Omega=\prod_{i=1}^{d} (\check{a}_i, \hat{a}_i)$; $\bm{x}=(x_1, x_2, \ldots, x_d) \in \mathbb{R}^d;$ for $i= 1, \ldots, d$ and $f,\psi$ are given functions. The operator \(\mathcal{L}\) represents the spatial differential operator acting on \(u(x,t)\), which can take several representative forms such as

\[
\mathcal{L} =
\begin{cases}
\Delta, & \text{constant Laplacian},\\[4pt]
\displaystyle
\sum_{i=1}^{d} c_i\,
\frac{\partial^{\alpha_i}}{\partial |x_i|^{\alpha_i}},
& \text{Riesz fractional derivative, } \alpha_i\in(1,2),\\[10pt]
\displaystyle
\sum_{i=1}^{d}
\left(
k_{i,+}\,
\frac{\partial^{\alpha_i}}{\partial_{+}x_i^{\alpha_i}}
+
k_{i,-}\,
\frac{\partial^{\alpha_i}}{\partial_{-}x_i^{\alpha_i}}
\right),
& \text{Riemann--Liouville fractional derivative, } \alpha_i\in(1,2).
\end{cases}
\]

Numerical methods for evolutionary FDEs of the form~\eqref{1d} with these choices of \(\mathcal{L}\) have been extensively discussed in the literature; see \cite{Chen2009,Gan2024,Jin2016,Lin2024,Lin2018,Lin2019,Pang2016,Vong2016,Zhu2024} and the references therein.\\
Because an explicit solution of~\eqref{1d} is rarely available, one usually relies on numerical discretization techniques; see, for example, \cite{Arshad2017,Bu2015,Chen2013,Feng2016,Hamid2019,Song2014,Sun2016,Yu2012,Yu2013,Zhu2018}. The fractional differential operators are nonlocal, and therefore the resulting discretized systems are dense, when using standard approximation methods. As a result, direct methods such as Gaussian elimination require high computational cost and become inefficient for large problems. However, the shift-invariant nature of the kernels in fractional operators leads to matrices with Toeplitz-type structure. This structure allows fast matrix-vector products through the use of fast Fourier transforms, with $O(n\log n)$ arithmetic cost, $n$ being the matrix size; see \cite{CN-sirev,VanLoan} and references therein. Using these efficient algorithms, iterative methods are designed to solve the resulting linear systems effectively, since at every step few matrix-vector products are performed.

The numerical methods for the resulting linear systems can generally be divided into two categories: time-stepping approaches and all-at-once approaches. In time-stepping methods, the system is solved sequentially in time. More precisely, the solution at the $n$-th time level must be obtained before computing the solution at the $(n+1)$-th time level. For representative works on time-stepping strategies, see \cite{Baffet2017,Jiang2017,Lin2019}. In contrast, all-at-once methods treat the entire space-time system simultaneously. These approaches compute the unknowns at all time levels together, which allows parallel implementation. For related studies on all-at-once solvers, see \cite{Gu2020,Lin2018,Zhao2019,Zhao2021}.

The present work focuses on the case $d=1$, and more precisely in the design of an efficient preconditioner for the all-at-once linear system arising from the discretization of \eqref{1d} for $d=1$. Since the preconditioner must approximate the coefficient matrix in a spectral sense to be effective, we first perform an asymptotic spectral analysis of the involved matrix sequences. This analysis is carried out using the theory of generalized locally Toeplitz (GLT) matrix sequences, which represents a powerful framework for studying the spectral behavior of matrices arising from variable-coefficient differential operators.
Notable applications of this theoretical apparatus comprise matrix sequences stemming from differential problems approximated by finite differences methods \cite{Serra2003}, finite elements and isogeometric analysis methods \cite[Chapter 10]{garoni2017}, \cite{systems1}, finite volumes techniques \cite{Bertaccini2017}, discontinuous Galerkin methods \cite{systems2} or stemming by the numerical approximation applied to fractional problems \cite{DMS18,Ilyas2026}.

The GLT theory is particularly well-suited for our problem because the spatially varying diffusion coefficient $a(x)$ spoils the pure Toeplitz structure that would otherwise be present in the constant-coefficient case; see \cite{barbarino2020uni, barbarino2020multi, garoni2017, garoni2018} for a complete treatment of the theory and the tutorials \cite{Garoni2018,Garoni2019} for a gentle guide on the use of the theory in practical applications. By characterizing the GLT symbols of the spatial and temporal discretization matrices, we are able to construct a preconditioner that captures the essential spectral features of the full coefficient matrix.

Specifically, we propose a block lower triangular preconditioner $\mathcal{P}$ that mimics the structure of the all-at-once matrix $\mathcal{A}$ while simplifying its components for computational efficiency. The preconditioner is designed such that the preconditioned matrix $\mathcal{P}^{-1}\mathcal{A}$ differs from the identity by a low-rank correction, ensuring a weak clustering of the eigenvalues at $1$, for the associated matrix sequence. This spectral clustering is rigorously analyzed within the GLT framework, and the theoretical findings are validated through numerical experiments using the GMRES solver.

The proposed technique and the related analysis are quite versatile and can adapted to the general problem in $d$ space dimensions, also of non-Cartesian type. In that case a (reduced) $(d+1)$-level GLT structure arise; see \cite{garoni2018,Serra2006,reduced} where the external level is related to discretization with respect to the time variable and the internal $d$ levels are related to the approximation formulae in the $d$ spatial variables.

The paper is organized as follows. Section~\ref{GLT Tool} introduces the theoretical framework and provides a summary of the essential tools. Section~\ref{time-space} presents the discretization of the nonlocal evolutionary FDE using finite difference schemes in time and space, and derives the resulting all-at-once linear system, emphasizing the structural features of the matrices. In Section~\ref{d+1 sys}, the coefficient matrix sequence is studied via GLT tools and its symbol is derived. Section \ref{algo-prec} is devoted to algorithmic proposals, to preconditioning strategies, and to the same type of spectral analysis for the preconditioned matrix sequences, proving clustering results both in the singular value and eigenvalue sense. In Section \ref{numerics}, numerical experiments are conducted to assess the effectiveness of the proposed preconditioner in terms of iteration counts and CPU time. Concluding remarks are presented in Section~\ref{conclusion} together with a sketch of the main ideas for the $d$ dimensional space case and further open problems.


\section{Preliminaries}\label{GLT Tool}
We concisely introduce the basic tools employed in this work, which are discussed in detail in \cite{garoni2017}. For simplicity, we restrict the definitions and theorems to the unilevel scalar case, which is sufficient for the purposes of the current paper. However, all the concepts can be extended to the multilevel \cite{garoni2018} and block \cite{barbarino2020uni, barbarino2020multi} settings, that arise from more complex differential problems such as systems of partial differential equations \cite{systems1,systems2} or more precise approximations such the high order finite elements \cite{GSS} and isogeometric analysis with intermediate regularity \cite{Garoni2019}.
We observe that the approximation of the $d$-dimensional space equation ~\eqref{1d} leads to $(d+1)$-level GLT matrix sequence with a specific tensor structure which can be treated using only unilevel GLT structures (see Axiom \textbf{GLT4}).

Let us first establish some useful notation.
\begin{itemize}
    \item Given a square matrix $A_n$, we denote with $\lambda_j(A_n)$ and $\sigma_j(A_n)$ the $j$-th eigenvalue and singular value of $A_n$, respectively, as $j = 1, \ldots, n$.
    \item The spectral norm for a square matrix $A_n$ is denoted as $\|A_n\|$.
    \item Whenever we use terminology from measure theory, for instance ``measurable set'', ``measurable function'', ``a.e.'', we always refer to the Lebesgue measure in $\mathbb{R}^t$, denoted with $\mu_t$.
\end{itemize}
\subsection{Spectral tools}
Throughout the paper, we deal with matrix sequences, referring to any sequence of
the form $\{A_n\}_n$, where $A_n$ is a square matrix of size $d_n$ and $d_n \to \infty$ monotonically
as $n \to \infty$. To a given matrix sequence, it is often possible to associate a spectral or
singular value distribution, according to the following definition.

\begin{Definition}\label{sig and eig }
Let $\{A_n\}_n$ be a matrix sequence, where $A_n$ has size $d_n \times d_n$, and let $\psi : D \subset \mathbb{R}^t \to \mathbb{C}$ be a measurable function defined on a set $D$ with $0 < \mu_t(D) < \infty$. Denote with $C_c(\mathbb{K})$ the set of continuous complex-valued functions with bounded support on $\mathbb{K} \in \{\mathbb{R}, \mathbb{C}\}$.

\begin{itemize}
    \item $\{A_n\}_n$ has an \textit{(asymptotic) singular value distribution} described by $\psi$ if for any $F \in C_c(\mathbb{R})$
    \[
    \lim_{n \to \infty} \frac{1}{d_n} \sum_{i=1}^{d_n} F(\sigma_i(A_n)) = \frac{1}{\mu_t(D)} \int_D F(|\psi(x)|) \, dx,
    \]
    where $\sigma_i(A_n)$, $i = 1, \ldots, d_n$, are the singular values of $A_n$. We use the notation $\{A_n\}_n \sim_\sigma \psi$.

    \item $\{A_n\}_n$ has an \textit{(asymptotic) spectral or eigenvalue distribution} described by $\psi$ if for any $F \in C_c(\mathbb{C})$
    \[
    \lim_{n \to \infty} \frac{1}{d_n} \sum_{i=1}^{d_n} F(\lambda_i(A_n)) = \frac{1}{\mu_t(D)} \int_D F(\psi(x)) \, dx,
    \]
    where $\lambda_i(A_n)$, $i = 1, \ldots, d_n$, are the eigenvalues of $A_n$. We use the notation $\{A_n\}_n \sim_\lambda \psi$.
\end{itemize}
\end{Definition}

The informal meaning behind Definition~\ref{sig and eig } is the following: if $\psi$ is continuous and $n$ is large, then the eigenvalues or singular values of $A_n$ (suitably ordered and possibly except at most for $o(d_n)$ outliers) are approximated by a uniform sampling of $\psi$ or $|\psi|$ over its domain; see \cite[Remark 2]{Ilyas2025} and also \cite[Remark 2.7]{MRS} for the more general setting of matrix-valued symbols.

\begin{Remark}
The relation $\{A_n\}_n \sim_{\lambda} \psi$ and $\Lambda(A_n) \subseteq S$ for all $n$ imply that the range of $\psi$ is a subset of the closure $\overline{S}$ of $S$. In particular, $\{A_n\}_n \sim_{\lambda} \psi$ and $A_n$ positive definite for all $n$ imply that $\psi$ is nonnegative almost everywhere; see e.g. \cite{Barbarino2020} for its use in the spectral analysis of approximated FDEs and \cite{Ahmad2025, Ilyas2025} for the spectral study of matrix sequences stemming from geometric means.
\end{Remark}

The notion of clustering can be seen a special case of distribution. In what follows, we denote the $\varepsilon$-expansion of a set $S \subseteq \mathbb{C}$ as
\[
B(S, \varepsilon) := \bigcup_{z \in S} B(z, \varepsilon),
\]
where $B(z, \varepsilon) := \{w \in \mathbb{C} : |w - z| < \varepsilon\}$ is the complex disk with center $z$ and radius $\varepsilon > 0$.

\begin{Definition}\label{cluster and distri}
Let $\{A_n\}_n$ be a matrix sequence, with $A_n$ of size $d_n \times d_n$, and let $S \subseteq \mathbb{C}$ be nonempty and closed. The sequence $\{A_n\}_n$ is \textbf{strongly clustered} at $S$ in the sense of the eigenvalues if $\forall \varepsilon > 0$
\[
\#\{j \in \{1, \ldots, d_n\} : \lambda_j(A_n) \notin B(S, \varepsilon)\} = O(1), \quad n \to \infty.
\]
In other words, the number of eigenvalues of $A_n$ outside $B(S, \varepsilon)$ is bounded by a constant independent of $n$. Moreover, $\{A_n\}_n$ is \textbf{weakly clustered} at $S$ if $\forall \varepsilon > 0$
\[
\#\{j \in \{1, \ldots, d_n\} : \lambda_j(A_n) \notin B(S, \varepsilon)\} = o(d_n), \quad n \to \infty,
\]
meaning that the number of eigenvalues of $A_n$ outside $B(S, \varepsilon)$ is negligible with respect to the size of the matrix. A corresponding definition can be given for the singular values, with $S \subseteq \mathbb{R}^+$.
\end{Definition}

\begin{Remark}\label{distribution}
Let us clarify the relationship between the concepts of distribution and clustering. We recall that the essential range of a measurable function $g : D \subseteq \mathbb{R}^t \to \mathbb{C}$ is the set
\[
\{ z \in \mathbb{C} : \mu_t(\{ g \in B(z, \varepsilon) \}) > 0 \quad \forall \varepsilon > 0 \}.
\]
As reported in \cite[Theorem 4.2]{Golinskii2007}, it holds
\[
\{ A_n \}_n \sim_\lambda \psi \quad \Longrightarrow \quad \{ A_n \}_n \text{ is weakly clustered at the essential range of } \psi.
\]
Furthermore, if the essential range of $\psi$ consists of a single number $s \in \mathbb{C}$, then
\[
\{ A_n \}_n \sim_\lambda \psi \quad \Longleftrightarrow \quad \{ A_n \}_n \text{ is weakly clustered at } s \text{ in the sense of the eigenvalues}.
\]
The latter case is typically of interest in the context of preconditioning. Corresponding statements can be given for the singular values, with the obvious suitable changes.
\end{Remark}

\subsection{\texorpdfstring{The GLT $*$-algebra}{The GLT *-algebra}}

A GLT sequence is a matrix sequence belonging to the $*$-algebra generated by the three specific classes of matrix sequences: zero-distributed, Toeplitz and diagonal sampling matrix sequences. We define them in the following paragraphs and they can be seen as the basic building blocks of the GLT $*$-algebra.

Each GLT sequence is equipped with a measurable function $\kappa : [0,1]^d \times [-\pi, \pi]^d \to \mathbb{C}$, with $d \geq 1$, called the GLT symbol. When $d = 1$, we are in the context of unilevel GLT sequences, while the case $d > 1$ relates to the multilevel GLT setting. The symbol is essentially unique, in the sense that if $\kappa, \zeta$ are two symbols of the same GLT sequence, then $\kappa = \zeta$ a.e. We write $\{A_n\}_n \sim_{\mathrm{GLT}} \kappa$ to denote that $\{A_n\}_n$ is a GLT sequence with symbol $\kappa$.

\paragraph{Zero-distributed sequences}
A matrix sequence $\{Z_n\}_n$ is zero-distributed if it holds $\{Z_n\}_n \sim_\sigma 0$. The following Lemma provides a practical characterization.

\begin{lemma}\label{lem:zero_dist}
Let $\{Z_n\}_n$ be a matrix sequence, with $Z_n$ of size $d_n \times d_n$. Then, $\{Z_n\}_n \sim_\sigma 0$ if and only if for every $n \in \mathbb{N}$ it holds $Z_n = R_n + N_n$ with
\[
\lim_{n \to \infty} \frac{\operatorname{rank}(R_n)}{d_n} = \lim_{n \to \infty} \|N_n\| = 0.
\]
\end{lemma}
\paragraph{Toeplitz sequences}
For any $n \in \mathbb{N}$, a matrix of the form
\[
[a_{s-t}]_{s,t=1}^n \in \mathbb{C}^{n \times n},
\]
with coefficients $a_k \in \mathbb{C}$ for $k = 1 - n, \ldots, n - 1$, is called a (unilevel) Toeplitz matrix and it is characterized by constant elements along the diagonals. Given a complex-valued function $f \in L^1([-\pi, \pi])$, the $n$-th Toeplitz matrix associated with $f$ is defined as
\begin{equation}\label{def:toe}
T_n(f) := [\hat{f}_{s-t}]_{s,t=1}^n \in \mathbb{C}^{n \times n},
\end{equation}
where
\begin{equation}\label{def:toe-fou}
\hat{f}_k = \frac{1}{2\pi} \int_{[-\pi, \pi]} f(\theta)e^{-ik\theta} d\theta, \quad k \in \mathbb{Z}, \quad i^2 = -1,
\end{equation}
are the Fourier coefficients of $f$. The family $\{T_n(f)\}_{n \in \mathbb{N}}$ is the sequence of Toeplitz matrices associated with $f$, called the generating function of the sequence.

We recall that the Wiener class is the sub-algebra of the continuous $2\pi$-periodic functions on $[-\pi, \pi]$ and it corresponds to the set of functions whose Fourier series is absolutely convergent, i.e.,
\[
\left\{ f : [-\pi, \pi] \to \mathbb{C} : \sum_{k \in \mathbb{Z}} |\hat{f}_k| < +\infty \right\}.
\]

\paragraph{Diagonal sampling sequences}
Given $a : [0, 1] \to \mathbb{C}$ continuous a.e., for any $n \in \mathbb{N}$ we define the $n$-th diagonal sampling matrix $D_n(a)$ as
\begin{equation}
D_n(a) = \operatorname{diag}_{i=1,\ldots,n} a\left( \frac{i}{n} \right) \in \mathbb{C}^{n \times n} \label{eq:dna}.
\end{equation}
The family $\{D_n(a)\}_{n \in \mathbb{N}}$ is the sequence of diagonal sampling matrices generated by $a(x)$.

The GLT class satisfies several algebraic and topological properties that are treated in great detail and generality in \cite{barbarino2020uni, barbarino2020multi, garoni2017, garoni2018}. For the purposes of this work, it is sufficient to present them through their operative properties, listed below in the unilevel setting i.e. for $d = 1$.

\vspace{0.4cm}
\textbf{GLT1.} If $\{A_n\}_n \sim_{\mathrm{GLT}} \kappa$, then $\{A_n\}_n \sim_\sigma \kappa$ in the sense of Definition~\ref{sig and eig }. If moreover each $A_n$ is Hermitian, then $\{A_n\}_n \sim_\lambda \kappa$.

\vspace{0.2cm}
\textbf{GLT2.} It holds

\begin{itemize}
    \item $\{T_n(f)\}_n \sim_{\mathrm{GLT}} \kappa(x, \theta) = f(\theta)$ for any $f \in L^1([-\pi, \pi])$;
    \item $\{D_n(a)\}_n \sim_{\mathrm{GLT}} \kappa(x, \theta) = a(x)$ for any $a : [0, 1] \to \mathbb{C}$ continuous a.e.;
    \item $\{Z_n\}_n \sim_{\mathrm{GLT}} \kappa(x, \theta) = 0$ if and only if $\{Z_n\}_n \sim_\sigma 0$, i.e., $\{Z_n\}_n$ is zero-distributed.
\end{itemize}

\textbf{GLT3.} If $\{A_n\}_n \sim_{\mathrm{GLT}} \kappa$ and $\{B_n\}_n \sim_{\mathrm{GLT}} \zeta$, then
\begin{itemize}
    \item $\{A_n^*\}_n \sim_{\mathrm{GLT}} \bar{\kappa}$, where $A_n^*$ denotes the conjugate transpose of $A_n$ and $\bar{\kappa}$ the complex conjugate of $\kappa$;
    \item $\{\alpha A_n + \beta B_n\}_n \sim_{\mathrm{GLT}} \alpha \kappa + \beta \zeta$ for any $\alpha, \beta \in \mathbb{C}$;
    \item $\{A_n B_n\}_n \sim_{\mathrm{GLT}} \kappa \zeta$;
    \item $\{A_n^\dagger\}_n \sim_{\mathrm{GLT}} \kappa^{-1}$ for any $\kappa$ invertible a.e., where $A_n^\dagger$ denotes the Moore-Penrose pseudoinverse of $A_n$.
\end{itemize}

\textbf{GLT4.} Let $\{A_n\}_n \sim_{\mathrm{GLT}} \kappa$ and $\{B_m\}_m \sim_{\mathrm{GLT}} \xi$ with
\[
\kappa(x_1, \theta_1) : [0, 1] \times [-\pi, \pi] \to \mathbb{C},
\]
\[
\xi(x_2, \theta_2) : [0, 1] \times [-\pi, \pi] \to \mathbb{C}.
\]
Then, setting $N = N(n, m) := nm$, it holds $\{A_n \otimes B_m\}_N \sim_{\mathrm{GLT}} \kappa \otimes \xi$, where $(\kappa \otimes \xi) : [0, 1]^2 \times [-\pi, \pi]^2 \to \mathbb{C}$ is given by
\[
(\kappa \otimes \xi)(x_1, x_2, \theta_1, \theta_2) := \kappa(x_1, \theta_1) \xi(x_2, \theta_2),
\]
i.e. we obtain a two-level GLT matrix sequence ($d = 2$). Notice that the axioms in the $d$-level setting are verbatim the same with $[0, 1]$ and $[-\pi, \pi]$ replaced by $[0, 1]^d$ and $[-\pi, \pi]^d$, respectively.
%
%

\section{\texorpdfstring{Discretization of~\eqref{1d}}{Discretization of (1d)} and the space-time linear system}
\label{time-space}

We divide the section into two parts, one devoted to the time discretization and one devoted to the space discretization.

\subsection{Temporal Discretizations}\label{temdis}
Let \(\mathbb{N}^{+}\) denote the set of all positive integers.
For \(M\in\mathbb{N}^{+}\), define the temporal step size \( \Delta t= {T}/{M}\) and \(t_m = m \Delta T\). Using the L1 scheme
(see, e.g., \cite{Jin2016,Liao2018,LinXu2007,SunWu2006}), the temporal discretization of the fractional derivative in~\eqref{1d} takes the form
\begin{align}\label{tempral}
\frac{1}{\Gamma(1-\beta)}
\int_{0}^{t}
\frac{\partial u({x},s)}{\partial s}\,
(t-s)^{-\beta}\,ds
\Big|_{\,t=t_m} \,\,\,
& = \
\frac{1}{\Gamma(2-\beta)} \frac{1}{(\Delta t)^\beta} \left [
 \ell_0^{(\beta)}\, u(x, t_m) + \sum_{j=1}^{m-1} \ell_j^{(\beta)}\,u({x},t_{m-j}) + \ell_m^{(\beta)}\,\psi({x}) \right ] \nonumber \\
& \quad +
O(\Delta t^{\,2-\beta}), \quad {x}\in \Omega,
\end{align}
where the coefficients are given by
\[
\ell_{j}^{(\beta)}=
\begin{cases}
 b_{0}, & j=0,\\[4pt]
(b_{j}-b_{j-1}), & j=1, \ldots, m-1,\\[4pt]
-b_{m-1}, & j=m.
\end{cases}
\]
%
\subsection{Spatial Discretization}\label{timedisc}
Suppose that a uniform spatial discretization with step size $\Delta x = ({\check{a}-\hat{a}})/({N+1})$ is adopted, with $N \in \mathbb{N}^{+}$. Under this discretization, we consider the following specific form of the operator $\mathcal{L} = a({x})\,{}_{a}D_{{x}}^{\alpha} u({x}, t),$ as introduced in \cite{khan2026structure}. The operator is approximated using the left sided Riemann–Liouville fractional derivative, discretized by means of the shifted Grünwald scheme, which leads to a nonsymmetric Toeplitz matrix having the following form
\[
 \overline{G}_{\alpha, N} =
\begin{bmatrix}
g_1^{(\alpha)} & g_0^{(\alpha)}             & 0             & \cdots & 0 \\
g_2^{(\alpha)} & g_1^{(\alpha)} & g_0^{(\alpha)}            & \ddots & \vdots \\
\vdots & g_2^{(\alpha)} & g_1^{(\alpha)} & \ddots & 0 \\
\vdots        & \ddots        & \ddots        & \ddots & g_0^{(\alpha)}  \\
g_{N}^{(\alpha)} &\cdots & \cdots & g_2^{(\alpha)} & g_1^{(\alpha)}
\end{bmatrix} \in \mathbb{R}^{N \times N}.
\]
The spatially varying diffusion coefficient $a({x})$ is incorporated through pointwise sampling
${\mathcal D}_{N}(a) = \mathrm{diag} \left( a ( x_i) \right)_{i=1,\ldots,N}$
%
In terms of notations, notice that ${\mathcal D}_{N}(a)$ coincides with $D_N(a)$ as in (\ref{eq:dna}) with $d=r=1$,
provided that the domain of $a$ is $[0,1]$. In general $[a_1,b_1]\neq [0,1]$ and hence ${\mathcal D}_{N}(a)=D_N(\hat a)$ with $\hat a(\hat x)=a(a_1 + (b_1-a_1)\hat x)$ and $x=a_1 + (b_1-a_1)\hat x$, $\hat x\in [0,1]$. This rescaling is useful for applying the GLT axioms when deducing the Weyl distributions of the related matrix sequences. 
Consequently, the overall spatial discretization matrix 
can be written as
\begin{equation}
    \frac{1}{(\Delta x)^{\alpha}} \mathcal{D}_{N}(a)\,\overline{G}_{\alpha, N},
\end{equation}
Since the construction of the discretizations and the analysis of the associated matrix properties have already been carried out in detail in \cite{khan2026structure}, we omit further technical details here and directly employ the resulting matrix formulation in the subsequent analysis.

\subsection{All At Once system}
By combining the spatial discretizations matrix 
with the temporal discretization scheme given in~\eqref{tempral}, for $t = t_m$ fixed and $x_i$, $i=1, \ldots, N$ we have
\[
\frac{1}{\Gamma(2-\beta)(\Delta t)^\beta}  \left [
 \ell_0^{(\beta)}\, \mathbf{u}^{(m)}_N + \sum_{j=1}^{m-1} \ell_j^{(\beta)}\,\mathbf{u}^{(m-j)}_N + \ell_m^{(\beta)}\,\boldsymbol{\psi}_N \right ]
 = \frac{1}{(\Delta x)^\alpha} \mathcal{D}_{N}(a)\,\overline{G}_{\alpha, N} \mathbf{u}^{(m)}_N +\mathbf{f}^{(m)}_N \\
\]
where
\( \mathbf{u}^{(m)}_N  = [u^{(m)}_i]^T_{i=1,\ldots,N}\) represents the numerical approximation of the solution $u(x,t_m)$
at the interior grid points $x_i$, $i=1,\ldots,N$, and accordingly the vectors \( \mathbf{f}^{(m)}_N  = [ f(x_i,t_m)]^T_{i=1,\ldots,N}\), and \(\boldsymbol{\psi}_N = [\psi(x_i)]^T_{i=1,\ldots,N}\).

The discretized form of problem~\eqref{1d} can be written as the following all at once linear system,
compactly expressed using Kronecker products in the following form
\begin{equation}\label{all_at_once1d}
  \mathcal{A}_{N,M} \, \bm{u}
:=
\left(
B_M  \otimes I_N -  I_M \otimes  \gamma \mathcal{D}_N(a) \overline{G}_{\alpha, N}
\right) \bm{u}
=
\bm{f},
\end{equation}
with \(\bm{u} = [\mathbf{u}^{(1)}_N, \ldots, \mathbf{u}^{(M)}_N]^T\) and
where the lower triangular Toeplitz matrix $B_M \in \mathbb{R}^{M \times M}$ denotes the temporal discretization matrix, namely,
\begin{equation}
B_M :=
\begin{bmatrix}
\ell_0^{(\beta)} \\
\ell_1^{(\beta)} & \ell_0^{(\beta)} \\
\vdots & \ddots & \ddots \\
\ell_{M-2}^{(\beta)} & \cdots & \ddots & \ell_0^{(\beta)} \\
\ell_{M-1}^{(\beta)} & \ell_{M-2}^{(\beta)} & \cdots & \ell_1^{(\beta)} & \ell_0^{(\beta)}
\end{bmatrix}
=
\begin{bmatrix}
b_0 \\
b_1 - b_0 & b_0 \\
\vdots & \ddots & \ddots \\
b_{M-2} - b_{M-3} & \cdots & \ddots & b_0 \\
b_{M-1} - b_{M-2} & b_{M-2} - b_{M-3} & \cdots & b_1 - b_0 & b_0
\end{bmatrix}.
\end{equation}
Furthemore, \(I_N\) and \(I_M\) denote the \(N\times N\) and \(M\times M\) identity matrices, respectively, and \[\gamma= \Gamma(2-\beta)(\Delta t)^\beta/(\Delta x)^\alpha.\]

\section{Analysis of the space-time matrix structures}\label{d+1 sys}
%

In the current section we present the spectral analysis of the coefficient matrix sequence using the GLT framework to describe the asymptotic behavior of the matrices. We first examine the sequence generated by the spatial block $\mathcal{D}_N(a) \overline{G}_{\alpha, N}$ then analyze the time block $B_{M}$ and finally study the complete coefficient matrix sequence.
\vspace{0.6cm}

%
%
\begin{theorem}[{\cite[Theorem 3.5]{khan2026structure}}]\label{th:1D_spatial_new}
If \( a : [0, 1] \to \mathbb{R} \) is continuous almost everywhere, then
\begin{itemize}
    \item[$a1)$]
    $\{ \mathcal{D}_N(a) \overline{G}_{\alpha, N} \}_N \sim_{\mathrm{GLT}} a(x)f_\alpha(\theta)$,
    \item[$a2)$]
    $\{ \mathcal{D}_N(a) \overline{G}_{\alpha, N} \}_N \sim_{\sigma} a(x)f_\alpha(\theta)$,
\end{itemize}
where
\[
f_\alpha(\theta)
=
e^{-\mathrm{i}\theta}\left(1 - e^{\mathrm{i}\theta}\right)^{\alpha},
\qquad \theta\in[-\pi,\pi].
\]
\end{theorem}
Since \(\mathcal{D}_N(a) \overline{G}_{\alpha, N}\) is not Hermitian, its eigenvalue distribution cannot be directly deduced from \textbf{GLT1}, even if additional tools are available; see, e.g.,~\cite{eig-nonH}. In the constant-coefficient case, the eigenvalue distribution satisfies
$\{\overline{G}_{\alpha, N}\}_n \sim_{\lambda} f_\alpha(\theta),$
as a consequence of quite sophisticate results on Toeplitz matrix sequences~\cite{eig-nonH}; see also Figure~3.1 in~\cite{khan2026structure}, where the eigenvalues closely match the samples of the symbol.

In the variable-coefficient case, numerical evidence (see Figure~3.2 in~\cite{khan2026structure}) shows a strong agreement between the eigenvalues and the samples of the symbol \(a(x)f_\alpha(\theta)\). However, a complete theoretical proof for the variable-coefficient case is not available and remains an open problem.
%


\begin{Remark}\label{rem:gamma_new}
Let \({G}_{\alpha, N} = \gamma \mathcal{D}_N(a) \overline{G}_{\alpha, N}\).
The inclusion of the scaling factor $\gamma
=
{\Gamma(2-\beta)\,(\Delta t)^\beta}/{(\Delta x)^\alpha}$
in the analysis is straightforward and extends the scope of Theorem~\ref{th:1D_spatial_new}. Assume that
\(a:[0,1]\to\mathbb{R}\) is continuous almost everywhere and that
\[
\lim_{\Delta t,\Delta x \to 0^+}\gamma(\Delta t,\Delta x)
=
\lim_{\Delta t,\Delta x\to 0^+}
\frac{\Gamma(2-\beta)\,(\Delta t)^\beta}{(\Delta x)^\alpha}
=
\gamma^*.
\]
Then
\begin{align}
\{G_N\}_N
&\sim_{\mathrm{GLT}}
\gamma^* a(x)f_\alpha(\theta),
\label{eq:GJ_GLT}
\\
\{G_N\}_N
&\sim_{\sigma}
\gamma^* a(x)f_\alpha(\theta).
\label{eq:GJ_sigma}
\end{align}
Alternatively, since we are interested in solving linear systems, it is often convenient to remove the scaling factor and consider the normalized sequence
\begin{align}
\{\gamma^{-1}G_N\}_N
&\sim_{\mathrm{GLT}}
a(x)f_\alpha(\theta),
\label{eq:norm_GJ_GLT}
\\
\{\gamma^{-1}G_N\}_N
&\sim_{\sigma}
a(x)f_\alpha(\theta).
\label{eq:norm_GJ_sigma}
\end{align}
These normalized relations require only that
\(a:[0,1]\to\mathbb{R}\) be continuous almost everywhere, which is equivalent to the Riemann integrability of \(a\).

Finally, statements \(a1)\), \(a2)\),
\eqref{eq:GJ_GLT}--\eqref{eq:GJ_sigma}, and
\eqref{eq:norm_GJ_GLT}--\eqref{eq:norm_GJ_sigma}
remain valid if \(a\) is complex-valued, although in applications \(a\) is typically assumed to be nonnegative, often strictly positive.
\end{Remark}
%
\begin{theorem}\label{th:1dtime}
There exists a continuous $2\pi$-periodic function $g_{\beta}(\theta)$ such that
\begin{equation*}
    \{B_M\}_{M} \sim_{\mathrm{GLT}, \sigma} g_{\beta}(\theta), \qquad \theta \in [-\pi,\pi]
\end{equation*}
Moreover, $\{B_M\}_{M} \sim_{\lambda}  b_0$.
\end{theorem}
\begin{proof}
 To demonstrate that $\{B_M\}_M$ is a GLT sequence, we claim that $B_M$ is a Toeplitz matrix generated by a Lebesgue integrable function function, i.e., there exists $g_\beta$ such that $B_M = T_M(g_\beta)$ for all $M \in \mathbb{N}$ according to the notation in (\ref{def:toe})-(\ref{def:toe-fou}). The idea is simple: we prove that the series
\[
g_\beta(\theta) := \sum_{k=0}^{\infty} \ell_k e^{ik\theta},
\]
where the coefficients $\ell_k$ are given in Section~\eqref{temdis}, belongs to the Wiener class, so that $g_\beta$ is well defined, continuous and $2\pi$-periodic, and the corresponding Toeplitz matrix sequence is precisely $\{B_M\}_M$ by definition. From the definition of the L1 coefficients in Section~\ref{temdis}, we have $\ell_k = (b_k - b_{k-1})$ for $k \ge 1$. Substituting $b_m = (j+1)^{1-\beta} - j^{1-\beta}$ and simplifying directly yields
\[
\ell_k =  \left[ (k+1)^{1-\beta} - 2k^{1-\beta} + (k-1)^{1-\beta} \right].
\]
Let us consider the auxiliary function $\psi_\beta(t) := t^{1-\beta}$, where  $\Delta t = T/M$. Using $k = t_k/\Delta t$ with $t_k = k\Delta t$, we obtain
\[
\ell_k =  (\Delta t)^{\beta-1} \left[ (t_k + \Delta t)^{1-\beta} - 2t_k^{1-\beta} + (t_k - \Delta t)^{1-\beta} \right].
\]
Then replace $(t_k + \Delta t)^{1-\beta}$ with $\psi_\beta(t_{k+1})$, $t_k^{1-\beta}$ with $\psi_\beta(t_k)$, and $(t_k - \Delta t)^{1-\beta}$ with $\psi_\beta(t_{k-1})$, giving
\[
\ell_k =  (\Delta t)^{\beta-1} \left[ \psi_\beta(t_{k+1}) - 2\psi_\beta(t_k) + \psi_\beta(t_{k-1}) \right].
\]
Dividing and multiplying by $(\Delta t)^2$,
\[
\ell_k =  (\Delta t)^{\beta+1} \left[ \frac{\psi_\beta(t_{k+1}) - 2\psi_\beta(t_k) + \psi_\beta(t_{k-1})}{(\Delta t)^2} \right].
\]
The fraction is a second-order central difference approximation of $\psi_\beta''(t_k)$. Hence,
\[
\ell_k =  (\Delta t)^{\beta+1} \left[ \psi_\beta''(t_k) + O((\Delta t)^2) \right].
\]
Since $\psi_\beta''(t) = \beta(\beta-1) t^{-1-\beta}$ and $t_k = k\Delta t$, we get
\[
\begin{aligned}
\ell_k &=   (\Delta t)^{\beta+1} \beta(\beta-1) \frac{1}{(k\Delta t)^{1+\beta}} + (\Delta t)^{\beta+1} O((\Delta t)^{2}) \\
&=  (\Delta t)^{\beta+1} \beta(\beta-1) \frac{1}{k^{1+\beta} (\Delta t)^{1+\beta}} + O((\Delta t)^{\beta+3}) \\
&= \kappa \beta(\beta-1) \frac{1}{k^{1+\beta}} + O((\Delta t)^{\beta+3}).
\end{aligned}
\]
Thus $\ell_k = \dfrac{ \beta(\beta-1)}{k^{1+\beta}} + O((\Delta t)^{\beta+3})$, so $\ell_k \sim \dfrac{C}{k^{1+\beta}}$ as $k \to \infty$ with $C =  \beta(\beta-1)$. Since $\beta \in (0,1)$, we have $1+\beta > 1$, and therefore the series
\[
\sum_{k=0}^{\infty} |\ell_k|
\]
is convergent. Consequently, $g_\beta$ belongs to the Wiener class. Hence $g_\beta$ is continuous and $2\pi$-periodic, and $B_M = T_M(g_\beta)$ by definition.

By \textbf{GLT1}, we conclude
\[
\{B_M\}_M \sim_{\text{GLT},\sigma} g_\beta(\theta).
\]
Finally, all the eigenvalues of $B_M$ are trivially equal to $\ell_0 =  b_0$ because of the triangular structure, we have
\[
\{B_M\}_M \sim_\lambda  b_0.
\]
This follows immediately from the definition of distribution.
\end{proof}
\begin{Remark}
    As it is evident from the final part of the proofs of Theorem~\eqref{th:1D_spatial_new} and Theorem~\ref{th:1dtime}, the fact that the involved generating functions belong to the Wiener class is not a goal, but it is instrumental in applying the GLT axioms for deducing in a fast way
the distribution in the sense of the singular values
\end{Remark}
We are now ready to compute the GLT symbol of the sequence associated with the full coefficient matrix, shown in \eqref{all_at_once1d}. To provide meaningful results, we introduce a mild assumption regarding the rate at which the discretization parameters tend to zero with respect to one another.

\begin{theorem}\label{coefmatrixglt}
Let $a:[0,1]\to\mathbb{R}$ be continuous almost everywhere and let $\gamma=\frac{\Gamma(2-\beta)(\Delta t)^\beta}{(\Delta x)^\alpha},$
$\gamma^*$ as in Remark~\ref{rem:gamma_new}. Then

\begin{enumerate}
\item[(c1)] $\{\mathcal{A}_{N,M}\}_{N,M}\sim_{\mathrm{GLT}}\gamma^* a(x)f_\alpha(\theta_1) + g_\beta(\theta_2),$

\item[(c2)] $\{\mathcal{A}_{N,M}\}_{N,M}\sim_{\sigma}\gamma^* a(x)f_\alpha(\theta_1) + g_\beta(\theta_2),$
\end{enumerate}

with
$\mathcal{A}_{N,M} =  B_M  \otimes I_N +  I_M \otimes \gamma \mathcal{D}_N(a) \overline{G}_{\alpha, N}$,
 where $f_\alpha(\theta)=e^{-i\theta}(1-e^{i\theta})^\alpha$ and $g_\beta$ is the symbol of $B_N$ from Theorem~\ref{th:1dtime}.
\end{theorem}

\begin{proof}
It is enough to prove $(c1)$, since $(c2)$ can be directly inferred from $(c1)$ and \textbf{GLT1}. From Theorem~\eqref{th:1D_spatial_new}, we already know that the sequence of matrices $\{G_N\}_N$ is distributed in the GLT sense as $\gamma^* a(x)f_\alpha(\theta_1)$. Moreover, by Theorem~\eqref{th:1dtime}, the sequence $\{B_M\}_M$ admits the GLT symbol $g_\beta(\theta_2)$. Therefore, by combining these two GLT symbols through the Kronecker sum structure of $\mathcal{A}_{N,M}$ (see Axiom \textbf{GLT4}) and invoking the $*$-algebra property of the GLT class Axiom \textbf{GLT3}, part 2,  we deduce that
\[
\{\mathcal{A}_{N,M}\}_{N,M}  \sim_{\mathrm{GLT}} \gamma^* a(x)f_\alpha(\theta_1) + g_\beta(\theta_2),
\]
which establishes $(c1)$. The statement $(c2)$ then follows directly from $(c1)$ and Axiom \textbf{GLT1}.
\end{proof}
\section{Algorithmic proposals and complexity analysis}\label{algo-prec}

Due to the lower triangular block structure of $\mathcal{A}_{N,M}$ two different solution approaches can be considered:
\begin{itemize}
  \item[\bf Technique 1:] the forward block substitution solution approach,
  \item[\bf Technique 2:] the global solution approach.
\end{itemize}
In the forward block substitution solution approach we have to solve, for $k = 1, \ldots, M$
\[(\ell_0 I_N - \gamma D_N(a) \overline{G}_{\alpha, N}) \mathbf{u}_N^{(k)}  =  \mathbf{f}_N^{(k)} - \sum_{i=1}^{k-1} \ell_k \mathbf{u}_N^{(i)},\]
where each forward step requires the solution of a linear system with the very same matrix
\[
A_N = \ell_0 I_N - \gamma \mathcal{D}_N(a) \overline{G}_{\alpha, N},
\]
with $\ell_0 = b_0 = 1$ for all $\beta$.
\ \\
This can be done by considering a preconditioned GMRES where the preconditioner is chosen as
\[P_{\alpha,N} =\ell_0 I_N - \gamma d\, \widehat{C}_{\alpha,N},\]
where
\[d = \frac{1}{N} \sum_{i=1}^{N} d_i,
\]
with $d_i$ are the diagonal entries of $\mathcal{D}_N(a)$,
and where
$\widehat{C}_{\alpha,N}$ is the circulant matrix generated by the vector $\left[g_1^{(\alpha)}, \ldots, g_{N-1}^{(\alpha)}, g_0^{(\alpha)}\right]^T$,
that is by the first column of $\overline{G}_{\alpha,N}$ where the last entry $g_N^{(\alpha)}$ has been replaced by the coefficient $g_0^{(\alpha)}$, which is nonnegigible while $g_N^{(\alpha)}$ tends to zero as $N$ tends to infinity. 
\\%
Otherwise, for $\omega>0$ but small, the preconditioner is chosen as
\[P_{\omega,\alpha,N} = \ell_0 I_N - \gamma d\, \widehat{C}_{\omega,\alpha,N}\]
 with $\widehat{C}_{\omega,\alpha,N}$ corresponding $\omega-$circulant matrix \cite{ng-book,koro} related to $\overline{G}_{\alpha, N}$ and hence to $\widehat{C}_{\alpha,N}$: to be precise the first column of $\widehat{C}_{\omega,\alpha,N}$ is exactly $\left[g_1^{(\alpha)}, \ldots, g_{N-1}^{(\alpha)}, g_0^{(\alpha)}/\omega\right]^T$ so that the entries on the first upper diagonal of  $\widehat{C}_{\omega,\alpha,N}$ equal those of the original matrix $\overline{G}_{\alpha, N}$ i.e.
 \[
\left(\widehat{C}_{\omega,\alpha,N}\right)_{j,j+1}= \left(\overline{G}_{\alpha, N}\right)_{j,j+1}=g_0^{(\alpha)}, \ \ \ j=1,\ldots, N-1. 
\]
The principal part of the computational cost is given by $M$ applications of preconditioned GMRES, where the preconditioner solution cost via FFT is given by $O(N\log N)$ according to \cite{CN-sirev}, plus $\sum_{i=1}^{M}(i-1)N$ matrix-vector products with diagonal matrices. To sum up the computational cost is  $\frac{1}{2}   NM^2+\sharp it\, M \ O(N\log N)$, where $\sharp it$ is the average number of preconditioned GMRES iterations required to solve the M systems with matrix $A_N$.
\par 
On the other hand, a global solution approach has been considered as well, where the linear system with matrix $\mathcal{A}_{N,M}$ is solved by a preconditioned GMRES, where the preconditioner is chosen in a similar way as $P_{\alpha,\beta,NM} = C_{\alpha,\beta,NM}$ circulant matrix whose first column is given by
\[
[\ell_0,\ldots, \ell_{M-1}]^T \otimes [1, 0, \ldots,0]^T - \gamma [ \widehat{C}_{\omega,\alpha,N} [1, 0, \ldots,0]^T, 0, \ldots, 0]^T
\]
or $P_{\omega,\alpha,\beta,MN} = C_{\omega,\alpha,\beta,MN}$ corresponding $\omega-$circulant matrix of $C_{\alpha,\beta,NM}$. \\
The principal part of the related computational cost is given by direct applications of preconditioned GMRES, where the preconditioner solution cost via FFT is given by $\sharp \tilde{it}\, O(MN\log NM)$ and $\sharp \tilde{it}$ is the  number of required iterations.

In all the cases it is worth noticing that the cardinality of iterations is the critical parameter. When the quantity  $\sharp {it}$ is constant with respect to $N$ and the quantity  $\sharp \tilde{it}$ is constant with respect to $N$ and $M$ then we end up fast methods, with optimal cost of
$ O(MN\log NM)$ arithmetic operations for the second proposal.

Now these quantities $\sharp {it}$ and $\sharp \tilde{it}$ depend on a favourable spectral distribution of the corresponding matrix sequences.
In the subsequent part of the section, we study the singular value and eigenvalue distributions of the indicated preconditioned matrix sequences.

\subsection{GLT analysis of the preconditioning proposals}\label{gltprecon}

We study the quality of the approximation resulting from the replacement of the Toeplitz structures by the considered circulant and $\omega-$circulant counterparts.

\subsubsection{GLT analysis of the preconditioning for {\bf Technique 1}}

We first remind that $\overline{G}_{\alpha, N}=T_N(f_\alpha)$ with $f_\alpha$ belonging to the Wiener class; see e.g. \cite{DMS18}. Now by classical results on Toeplitz preconditoning \cite{CN-sirev,ng-book}, we know that the natural circulant and the natural $\omega$-circulant approximations $P_N$ ensure a strong singular value clustering at zero for 
$\{\overline{G}_{\alpha, N}-P_N\}_N$ in the sense of Definition \ref{cluster and distri}. Hence, taking into account Remark \ref{distribution}, $\{\overline{G}_{\alpha, N}-P_N\}_N\sim_{\sigma} 0$ and the latter is equivalent to
\begin{equation}\label{key-rel1}
\{\overline{G}_{\alpha, N}-P_N\}_N\sim_{\mathrm{GLT}} 0,
\end{equation}
by Axiom \textbf{GLT1}.

Relation (\ref{key-rel1}) is the basis for proving the following result.

\begin{theorem}\label{th1-prec}
Let $\alpha \in (1,2)$ and let $P_N\in \left\{\widehat{C}_{\alpha,N}, \widehat{C}_{\omega,\alpha,N}\right\}$. Then 
\begin{itemize}
\item $\{P_N\}_N\sim_{\mathrm{GLT},\sigma} f_\alpha(\theta)$;
\item $\{P_N^{-1}A_N\}_N\sim_{\mathrm{GLT},\sigma} 1$ when $a$ is a constant;
\item
\[
\{P_N^{-1}A_N\}_N \sim_{\mathrm{GLT},\sigma} \frac{\ell_0-\gamma a(x)f_\alpha(\theta)}{\ell_0-\gamma d f_\alpha(\theta)}
\]
for a general Riemann integrable coefficient $a(x)$.
\end{itemize}
\end{theorem}

\begin{proof}
By Axiom \textbf{GLT2}, part 1, we have $\{\overline{G}_{\alpha, N}\}_N\sim_{\mathrm{GLT}} f_\alpha(\theta)$ since $\overline{G}_{\alpha, N}=T_N(f_\alpha)$.
Then  $\{D_N(a)\overline{G}_{\alpha, N}\}_N\sim_{\mathrm{GLT}} a(x)f_\alpha(\theta)$ because of Axiom \textbf{GLT2}, part 2, and \textbf{GLT3}, part 3.
As a consequence, since $\{A_N\}_N$ is linear combination of the GLT sequences $\{I_N\}_N$ and $\{D_N(a)\overline{G}_{\alpha, N})\}_N$ with coefficients 
$\ell_0$ and $\gamma$, we deduce that
\begin{equation}\label{A_N GLT}
\{ A_N = \ell_0 I_N - \gamma \mathcal{D}_N(a) \overline{G}_{\alpha, N} \}_N\sim_{\mathrm{GLT}} \ell_0-\gamma a(x)f_\alpha(\theta).
\end{equation}
Now relation (\ref{key-rel1}), $\{\overline{G}_{\alpha, N}\}_N\sim_{\mathrm{GLT}} f_\alpha(\theta)$, and Axiom \textbf{GLT3}, part 2, imply 
 $\{P_N\}_N\sim_{\mathrm{GLT}} f_\alpha(\theta)$.
 Since $f_\alpha$ has a unique zero at $\theta=0$, it is obvious that it is invertible a.e. and hence Axiom \textbf{GLT3}, part 4, leads to the conclusion 
\[
\{P_N^{-1}\}_N\sim_{\mathrm{GLT}} \frac{1}{f_\alpha(\theta)}
\] 
so that, by Axiom \textbf{GLT3}, part 3, we deduce
\[
\{P_N^{-1}A_N\}_N\sim_{\mathrm{GLT}}=1
\]
if $a(x)$ is a constant function and 
\[
\{P_N^{-1}A_N\}_N \sim_{\mathrm{GLT}} \frac{\ell_0-\gamma a(x)f_\alpha(\theta)}{\ell_0-\gamma d f_\alpha(\theta)}
\]
when $a(x)$ is Riemann integrable.
Finally the corresponding singular value distributions follow directly from Axiom \textbf{GLT1}.
\end{proof}

\subsubsection{GLT analysis of the preconditioning for {\bf Technique 2}}

By Theorem \ref{th:1dtime}, we have $B_{M}=T_N(g_\beta)$ with $g_\beta$ belonging to the Wiener class. Following the literature on Toeplitz preconditoning \cite{CN-sirev,ng-book}, we know that the natural circulant and the natural $\omega$-circulant approximations $P_M$ ensure a strong singular value clustering at zero for
$\{B_M-P_M\}_M$ in the sense of Definition \ref{cluster and distri}. Hence, taking into account Remark \ref{distribution}, $\{B_M-P_M\}_M\sim_{\sigma} 0$ and the latter is equivalent to
\begin{equation}\label{key-rel2}
\{B_M-P_M\}_M\sim_{\mathrm{GLT}} 0,
\end{equation}
by Axiom \textbf{GLT1}.

Relation (\ref{key-rel2}) is the basis for proving the following result.

\begin{theorem}\label{th2-prec}
Let $\beta \in (1,2)$ and let $P_M\in \left\{\widehat{C}_{\beta,M}, \widehat{C}_{\omega,\beta,M}\right\}$. Then
\begin{itemize}
\item $\{P_M\}_M\sim_{\mathrm{GLT},\sigma} g_\beta(\theta_2)$;
\item $\{P_M^{-1}B_N\}_N\sim_{\mathrm{GLT},\sigma} 1$.
\end{itemize}
\end{theorem}
\begin{proof}
The proof is a simplified version of that of Theorem \ref{th1-prec}, which has been given with all the details, and hence we omit it.
\end{proof}

Joining the information in the last two theorem we can prove the following result, which is essentially based on Axiom \textbf{GLT4} on tensor product matrix sequences.

\begin{theorem}\label{th3-prec}
Let $\alpha,\beta \in (1,2)$. Then
\begin{itemize} 
\item $\{P_{\alpha,\beta,NM}\}_{N,M}\sim_{\mathrm{GLT},\sigma} \gamma^* d f_\alpha(\theta_1)+ g_\beta(\theta_2)$;
\item $\{P_{\alpha,\beta,NM}^{-1}\mathcal{A}_{N,M} \}_{N,M}\sim_{\mathrm{GLT},\sigma} 1$ when $a$ is a constant;
\item
\[
\{P_{\alpha,\beta,NM}^{-1}\mathcal{A}_{N,M} \}_{N,M}\sim_{\mathrm{GLT},\sigma} \frac{\gamma^* d a(x)f_\alpha(\theta_1)+ g_\beta(\theta_2)}{\gamma^* d f_\alpha(\theta_1)+ g_\beta(\theta_2)}
\]
for a general Riemann integrable coefficient $a(x)$.
\end{itemize}
\end{theorem}

\begin{proof}
The proof amounts in combining Theorem \ref{th1-prec}, Theorem \ref{th2-prec} with the tensor product Axiom \textbf{GLT4}.
\end{proof}


\section{Numerical Results}\label{numerics}

In the present section we perform extensive numerical experiments concerning the Krylov precondition related to {\bf Technique 1} and {\bf Technique 2}.
More precisely, Table \ref{tab:iterations_es1_alpha11}--Table \ref{tab:iterations_es2_alpha19} concern the number of the preconditioned GMRES in {\bf Technique 1} and Table \ref{tab:iterations_tot_es1_alpha11}--Table \ref{tab:iterations_tot_es2_alpha19} concern the number of the preconditioned GMRES in {\bf Technique 2} for the global linear systems.

In all the considered cases, irrespectively of the chosen parameters $\alpha$, $\beta$, $a(x)$, we can appreciated that the number of iterations stabilizes to a given constant independent of the matrix size. These constants are mild and are a bit more pronounced when $\beta$ is close to $2$.
The optimality of the preconditioning preconditioned GMRES is a confirmation that the total arithmetic cost is optimal with {\bf Technique 2}, since 
it grows in order as a bivariate FFT.


\begin{table}[H]
\centering
\vspace{0.5em}
\small
\begin{tabular}{|cc|ccc|ccc|ccc|}
\hline
& & \multicolumn{9}{c|}{$\alpha=1.1$} \\
& & \multicolumn{3}{c|}{$\beta=0.1$} & \multicolumn{3}{c|}{$\beta=0.5$} & \multicolumn{3}{c|}{$\beta=0.9$} \\
\cline{3-11}
$N$ & $M$
& -  & ${P_{\alpha,N}}$ & $P_{\omega,\alpha,N}$
& -  & ${P_{\alpha,N}}$ & $P_{\omega,\alpha,N}$
& -  & ${P_{\alpha,N}}$ & $P_{\omega,\alpha,N}$
\\
\hline
$2^4$ & $2^4$ & 16.00 & 4.00 & 4.00 & 16.00 & 4.00 & 4.00 & 16.00 & 4.00 & 4.00 \\
$2^5$ & $2^4$ & 32.00 & 4.00 & 4.00 & 32.00 & 5.00 & 4.00 & 32.00 & 4.00 & 4.00 \\
$2^6$ & $2^4$ & 64.00 & 4.00 & 4.00 & 64.00 & 5.00 & 4.00 & 64.00 & 4.06 & 4.00 \\
$2^7$ & $2^4$ & 128.00 & 5.00 & 4.00 & 128.00 & 5.00 & 4.00 & 128.00 & 5.00 & 4.00 \\
$2^8$ & $2^4$ & 256.00 & 5.00 & 4.00 & 256.00 & 5.00 & 4.00 & 256.00 & 5.00 & 4.00 \\
\hline
$2^4$ & $2^5$ & 16.00 & 4.00 & 4.00 & 16.00 & 4.00 & 4.00 & 16.00 & 4.00 & 4.00 \\
$2^5$ & $2^5$ & 32.00 & 4.00 & 4.00 & 32.00 & 4.09 & 4.00 & 32.00 & 4.00 & 4.00 \\
$2^6$ & $2^5$ & 64.00 & 4.00 & 4.00 & 64.00 & 5.00 & 4.00 & 62.12 & 4.00 & 4.00 \\
$2^7$ & $2^5$ & 128.00 & 5.00 & 4.00 & 128.00 & 5.00 & 4.00 & 117.22 & 4.09 & 4.00 \\
$2^8$ & $2^5$ & 256.00 & 5.00 & 4.00 & 256.00 & 5.00 & 4.00 & 225.31 & 5.00 & 4.00 \\
\hline
$2^4$ & $2^6$ & 16.00 & 4.00 & 4.00 & 16.00 & 4.00 & 4.00 & 14.00 & 4.00 & 4.00 \\
$2^5$ & $2^6$ & 32.00 & 4.00 & 4.00 & 32.00 & 4.02 & 4.00 & 22.03 & 4.00 & 4.00 \\
$2^6$ & $2^6$ & 64.00 & 5.00 & 4.00 & 64.00 & 5.00 & 4.00 & 38.05 & 4.00 & 4.00 \\
$2^7$ & $2^6$ & 128.00 & 5.00 & 4.00 & 128.00 & 5.00 & 4.00 & 70.03 & 4.00 & 4.00 \\
$2^8$ & $2^6$ & 256.00 & 5.00 & 4.00 & 256.00 & 5.00 & 4.00 & 132.06 & 4.03 & 4.00 \\
\hline
$2^4$ & $2^7$ & 16.00 & 4.00 & 4.00 & 16.00 & 4.00 & 4.00 & 10.00 & 4.00 & 4.00 \\
$2^5$ & $2^7$ & 32.00 & 4.00 & 4.00 & 32.00 & 4.00 & 4.00 & 15.00 & 4.00 & 4.00 \\
$2^6$ & $2^7$ & 64.00 & 5.00 & 4.00 & 64.00 & 4.01 & 4.00 & 24.01 & 4.00 & 4.00 \\
$2^7$ & $2^7$ & 128.00 & 5.00 & 4.00 & 128.00 & 5.00 & 4.00 & 42.01 & 4.00 & 4.00 \\
$2^8$ & $2^7$ & 256.00 & 5.00 & 4.00 & 256.00 & 5.00 & 4.00 & 77.03 & 4.00 & 4.00 \\
\hline
$2^4$ & $2^8$ & 16.00 & 4.00 & 4.00 & 16.00 & 4.00 & 4.00 & 8.00 & 4.00 & 3.00 \\
$2^5$ & $2^8$ & 32.00 & 4.00 & 4.00 & 32.00 & 4.00 & 4.00 & 11.00 & 4.00 & 4.00 \\
$2^6$ & $2^8$ & 64.00 & 5.00 & 4.00 & 64.00 & 4.00 & 4.00 & 16.00 & 4.00 & 4.00 \\
$2^7$ & $2^8$ & 128.00 & 5.00 & 4.00 & 128.00 & 5.00 & 4.00 & 26.00 & 4.00 & 4.00 \\
$2^8$ & $2^8$ & 256.00 & 5.00 & 4.00 & 256.00 & 5.00 & 4.00 & 46.00 & 4.00 & 4.00 \\
\hline
\end{tabular}
\caption{Averaged number of GMRES iterations in forward block solution without preconditioner, with $P_{\alpha,N}$ circulant preconditioner, and $P_{\omega,\alpha,N}$  $\omega-$circulant preconditioner - case $a(x) = 1$.}
\label{tab:iterations_es1_alpha11}
\end{table}
\begin{table}[H]
\centering
\vspace{0.5em}
\small
\begin{tabular}{|cc|ccc|ccc|ccc|}
\hline
& & \multicolumn{9}{c|}{$\alpha=1.5$} \\
& & \multicolumn{3}{c|}{$\beta=0.1$} & \multicolumn{3}{c|}{$\beta=0.5$} & \multicolumn{3}{c|}{$\beta=0.9$} \\
\cline{3-11}
$N$ & $M$
& -  & ${P_{\alpha,N}}$ & $P_{\omega,\alpha,N}$
& -  & ${P_{\alpha,N}}$ & $P_{\omega,\alpha,N}$
& -  & ${P_{\alpha,N}}$ & $P_{\omega,\alpha,N}$
\\
\hline
$2^4$ & $2^4$ & 16.00 & 5.00 & 4.00 & 16.00 & 5.00 & 4.00 & 16.00 & 5.00 & 4.00 \\
$2^5$ & $2^4$ & 32.00 & 5.00 & 4.00 & 32.00 & 5.00 & 4.00 & 32.00 & 5.00 & 4.00 \\
$2^6$ & $2^4$ & 64.00 & 5.00 & 4.00 & 64.00 & 5.00 & 4.00 & 64.00 & 5.00 & 4.00 \\
$2^7$ & $2^4$ & 128.00 & 5.00 & 4.00 & 128.00 & 5.00 & 4.00 & 128.00 & 5.00 & 4.00 \\
$2^8$ & $2^4$ & 256.00 & 5.00 & 4.00 & 256.00 & 5.00 & 4.00 & 256.00 & 5.00 & 4.00 \\
\hline
$2^4$ & $2^5$ & 16.00 & 5.00 & 4.00 & 16.00 & 5.00 & 4.00 & 16.00 & 5.00 & 4.00 \\
$2^5$ & $2^5$ & 32.00 & 5.00 & 4.00 & 32.00 & 5.00 & 4.00 & 32.00 & 5.00 & 4.00 \\
$2^6$ & $2^5$ & 64.00 & 4.97 & 4.00 & 64.00 & 5.00 & 4.00 & 64.00 & 5.00 & 4.00 \\
$2^7$ & $2^5$ & 128.00 & 5.00 & 4.00 & 128.00 & 5.00 & 4.00 & 128.00 & 5.00 & 4.00 \\
$2^8$ & $2^5$ & 256.00 & 5.00 & 4.00 & 256.00 & 5.00 & 4.00 & 256.00 & 4.97 & 4.00 \\
\hline
$2^4$ & $2^6$ & 16.00 & 5.00 & 4.00 & 16.00 & 5.00 & 4.00 & 16.00 & 5.00 & 4.00 \\
$2^5$ & $2^6$ & 32.00 & 5.00 & 4.00 & 32.00 & 5.00 & 4.00 & 32.00 & 5.00 & 4.00 \\
$2^6$ & $2^6$ & 64.00 & 4.98 & 4.00 & 64.00 & 5.00 & 4.00 & 64.00 & 5.00 & 4.00 \\
$2^7$ & $2^6$ & 128.00 & 4.98 & 4.00 & 128.00 & 5.00 & 4.00 & 128.00 & 5.00 & 4.00 \\
$2^8$ & $2^6$ & 256.00 & 4.97 & 4.00 & 256.00 & 5.00 & 4.00 & 256.00 & 5.00 & 4.00 \\
\hline
$2^4$ & $2^7$ & 16.00 & 5.00 & 4.00 & 16.00 & 5.00 & 4.00 & 16.00 & 4.00 & 4.00 \\
$2^5$ & $2^7$ & 32.00 & 5.00 & 4.00 & 32.00 & 5.00 & 4.00 & 30.01 & 5.00 & 4.00 \\
$2^6$ & $2^7$ & 64.00 & 4.99 & 4.00 & 64.00 & 5.00 & 4.00 & 56.02 & 5.00 & 4.00 \\
$2^7$ & $2^7$ & 128.00 & 5.00 & 4.00 & 128.00 & 5.00 & 4.00 & 107.01 & 5.00 & 4.00 \\
$2^8$ & $2^7$ & 256.00 & 4.98 & 4.00 & 256.00 & 5.00 & 4.00 & 203.07 & 5.00 & 4.00 \\
\hline
$2^4$ & $2^8$ & 16.00 & 5.00 & 4.00 & 16.00 & 5.00 & 4.00 & 12.01 & 4.00 & 4.00 \\
$2^5$ & $2^8$ & 32.00 & 5.00 & 4.00 & 32.00 & 5.00 & 4.00 & 21.00 & 4.00 & 4.00 \\
$2^6$ & $2^8$ & 64.00 & 5.00 & 4.00 & 64.00 & 5.00 & 4.00 & 38.05 & 5.00 & 4.00 \\
$2^7$ & $2^8$ & 128.00 & 5.00 & 4.00 & 128.00 & 5.00 & 4.00 & 73.01 & 5.00 & 4.00 \\
$2^8$ & $2^8$ & 256.00 & 4.99 & 4.00 & 256.00 & 5.00 & 4.00 & 139.03 & 5.00 & 4.00 \\
\hline
\end{tabular}
\caption{Averaged number of GMRES iterations in forward block solution without preconditioner, with $P_{\alpha,N}$ circulant preconditioner, and $P_{\omega,\alpha,N}$  $\omega-$circulant preconditioner - case $a(x) = 1$.}
\label{tab:iterations_es1_alpha15}
\end{table}
\begin{table}[H]
\centering
\vspace{0.5em}
\small
\begin{tabular}{|cc|ccc|ccc|ccc|}
\hline
& & \multicolumn{9}{c|}{$\alpha=1.9$} \\
& & \multicolumn{3}{c|}{$\beta=0.1$} & \multicolumn{3}{c|}{$\beta=0.5$} & \multicolumn{3}{c|}{$\beta=0.9$} \\
\cline{3-11}
$N$ & $M$
& -  & ${P_{\alpha,N}}$ & $P_{\omega,\alpha,N}$
& -  & ${P_{\alpha,N}}$ & $P_{\omega,\alpha,N}$
& -  & ${P_{\alpha,N}}$ & $P_{\omega,\alpha,N}$
\\
\hline
$2^4$ & $2^4$ & 16.00 & 4.00 & 4.00 & 16.00 & 4.00 & 4.00 & 16.00 & 4.00 & 4.00 \\
$2^5$ & $2^4$ & 32.00 & 4.00 & 4.00 & 32.00 & 4.00 & 4.00 & 32.00 & 5.00 & 4.00 \\
$2^6$ & $2^4$ & 64.00 & 4.00 & 4.00 & 64.00 & 4.00 & 4.00 & 64.00 & 5.00 & 4.00 \\
$2^7$ & $2^4$ & 128.00 & 4.00 & 4.00 & 128.00 & 4.00 & 4.00 & 128.00 & 5.00 & 4.00 \\
$2^8$ & $2^4$ & 256.00 & 4.00 & 4.00 & 256.00 & 4.00 & 4.00 & 256.00 & 4.25 & 4.00 \\
\hline
$2^4$ & $2^5$ & 16.00 & 4.00 & 4.00 & 16.00 & 4.00 & 4.00 & 16.00 & 4.00 & 4.00 \\
$2^5$ & $2^5$ & 32.00 & 4.00 & 4.00 & 32.00 & 4.03 & 4.00 & 32.00 & 5.00 & 4.00 \\
$2^6$ & $2^5$ & 64.00 & 4.00 & 4.00 & 64.00 & 4.06 & 4.00 & 64.00 & 5.00 & 4.00 \\
$2^7$ & $2^5$ & 128.00 & 4.00 & 4.00 & 128.00 & 4.03 & 4.00 & 128.00 & 5.00 & 4.00 \\
$2^8$ & $2^5$ & 256.00 & 4.00 & 4.00 & 256.00 & 4.00 & 4.00 & 256.00 & 5.00 & 4.00 \\
\hline
$2^4$ & $2^6$ & 16.00 & 4.00 & 4.00 & 16.00 & 4.00 & 4.00 & 16.00 & 4.00 & 4.00 \\
$2^5$ & $2^6$ & 32.00 & 4.00 & 4.00 & 32.00 & 4.19 & 4.00 & 32.00 & 5.00 & 4.00 \\
$2^6$ & $2^6$ & 64.00 & 4.00 & 4.00 & 64.00 & 5.00 & 4.00 & 64.00 & 5.00 & 4.00 \\
$2^7$ & $2^6$ & 128.00 & 4.00 & 4.00 & 128.00 & 4.61 & 4.00 & 128.00 & 5.00 & 4.00 \\
$2^8$ & $2^6$ & 256.00 & 4.00 & 4.00 & 256.00 & 4.02 & 4.00 & 256.00 & 5.00 & 4.00 \\
\hline
$2^4$ & $2^7$ & 16.00 & 4.00 & 4.00 & 16.00 & 4.00 & 4.00 & 16.00 & 3.99 & 4.00 \\
$2^5$ & $2^7$ & 32.00 & 4.00 & 4.00 & 32.00 & 5.00 & 4.00 & 32.00 & 4.02 & 4.00 \\
$2^6$ & $2^7$ & 64.00 & 4.00 & 4.00 & 64.00 & 5.00 & 4.00 & 64.00 & 5.00 & 4.00 \\
$2^7$ & $2^7$ & 128.00 & 4.00 & 4.00 & 128.00 & 5.00 & 4.00 & 128.00 & 5.00 & 4.00 \\
$2^8$ & $2^7$ & 256.00 & 4.00 & 4.00 & 256.00 & 4.03 & 4.00 & 256.00 & 5.00 & 4.00 \\
\hline
$2^4$ & $2^8$ & 16.00 & 3.97 & 4.00 & 16.00 & 4.00 & 4.00 & 16.00 & 4.00 & 4.00 \\
$2^5$ & $2^8$ & 32.00 & 4.00 & 4.00 & 32.00 & 4.99 & 4.00 & 32.00 & 4.00 & 4.00 \\
$2^6$ & $2^8$ & 64.00 & 4.00 & 4.00 & 64.00 & 4.96 & 4.00 & 64.00 & 4.02 & 4.00 \\
$2^7$ & $2^8$ & 128.00 & 4.00 & 4.00 & 128.00 & 5.00 & 4.00 & 125.05 & 4.93 & 4.00 \\
$2^8$ & $2^8$ & 256.00 & 4.00 & 4.00 & 256.00 & 5.00 & 4.00 & 236.30 & 4.01 & 4.00 \\
\hline
\end{tabular}
\caption{Averaged number of GMRES iterations in forward block solution without preconditioner, with $P_{\alpha,N}$ circulant preconditioner, and $P_{\omega,\alpha,N}$  $\omega-$circulant preconditioner - case $a(x) = 1$.}
\label{tab:iterations_es1_alpha19}
\end{table}
\begin{table}[H]
\centering
\vspace{0.5em}
\small
\begin{tabular}{|cc|ccc|ccc|ccc|}
\hline
& & \multicolumn{9}{c|}{$\alpha=1.1$} \\
& & \multicolumn{3}{c|}{$\beta=0.1$} & \multicolumn{3}{c|}{$\beta=0.5$} & \multicolumn{3}{c|}{$\beta=0.9$} \\
\cline{3-11}
$N$ & $M$
& -  & ${P_{\alpha,N}}$ & $P_{\omega,\alpha,N}$
& -  & ${P_{\alpha,N}}$ & $P_{\omega,\alpha,N}$
& -  & ${P_{\alpha,N}}$ & $P_{\omega,\alpha,N}$
\\
\hline
$2^4$ & $2^4$ & 16.00 & 7.81 & 8.06 & 16.00 & 8.00 & 7.06 & 16.00 & 7.06 & 8.00 \\
$2^5$ & $2^4$ & 32.00 & 8.00 & 9.00 & 32.00 & 8.00 & 8.06 & 32.00 & 8.12 & 9.00 \\
$2^6$ & $2^4$ & 64.00 & 8.00 & 9.06 & 64.00 & 8.88 & 8.06 & 64.00 & 9.00 & 9.00 \\
$2^7$ & $2^4$ & 128.00 & 8.00 & 9.06 & 128.00 & 9.00 & 8.94 & 128.00 & 9.00 & 9.00 \\
$2^8$ & $2^4$ & 256.00 & 8.00 & 9.06 & 256.00 & 9.00 & 8.94 & 256.00 & 9.06 & 9.06 \\
$2^4$ & $2^5$ & 16.00 & 7.91 & 8.03 & 16.00 & 8.00 & 8.00 & 16.00 & 7.00 & 7.00 \\
$2^5$ & $2^5$ & 32.00 & 8.00 & 9.00 & 32.00 & 8.03 & 8.03 & 32.00 & 8.00 & 8.97 \\
$2^6$ & $2^5$ & 64.00 & 8.00 & 9.03 & 64.00 & 9.00 & 9.00 & 64.00 & 9.00 & 9.00 \\
$2^7$ & $2^5$ & 128.00 & 8.00 & 9.03 & 128.00 & 9.00 & 9.00 & 128.00 & 9.03 & 9.97 \\
$2^8$ & $2^5$ & 256.00 & 8.00 & 9.03 & 256.00 & 9.00 & 9.00 & 256.00 & 9.06 & 10.00 \\
$2^4$ & $2^6$ & 16.00 & 7.97 & 8.02 & 16.00 & 7.98 & 8.00 & 16.00 & 6.98 & 7.00 \\
$2^5$ & $2^6$ & 32.00 & 8.00 & 9.00 & 32.00 & 8.03 & 8.97 & 27.02 & 8.00 & 8.00 \\
$2^6$ & $2^6$ & 64.00 & 8.00 & 9.02 & 64.00 & 9.00 & 9.00 & 49.02 & 9.00 & 9.00 \\
$2^7$ & $2^6$ & 128.00 & 8.00 & 9.02 & 128.00 & 9.00 & 9.00 & 92.02 & 9.02 & 9.98 \\
$2^8$ & $2^6$ & 256.00 & 8.00 & 9.02 & 256.00 & 9.00 & 9.00 & 176.05 & 9.98 & 10.00 \\
$2^4$ & $2^7$ & 16.00 & 7.98 & 8.01 & 16.00 & 7.01 & 8.00 & 12.00 & 6.00 & 6.00 \\
$2^5$ & $2^7$ & 32.00 & 8.00 & 9.00 & 32.00 & 8.02 & 9.00 & 19.01 & 7.00 & 8.00 \\
$2^6$ & $2^7$ & 64.00 & 8.00 & 9.01 & 64.00 & 9.00 & 9.00 & 33.01 & 8.98 & 9.00 \\
$2^7$ & $2^7$ & 128.00 & 8.00 & 9.01 & 128.00 & 9.00 & 9.00 & 60.02 & 9.00 & 9.99 \\
$2^8$ & $2^7$ & 256.00 & 8.00 & 9.01 & 256.00 & 9.01 & 9.01 & 114.04 & 10.00 & 10.00 \\
$2^4$ & $2^8$ & 16.00 & 7.99 & 8.00 & 16.00 & 7.00 & 8.00 & 9.00 & 5.98 & 6.00 \\
$2^5$ & $2^8$ & 32.00 & 8.00 & 9.00 & 32.00 & 8.00 & 9.00 & 14.00 & 7.00 & 7.00 \\
$2^6$ & $2^8$ & 64.00 & 8.00 & 9.00 & 64.00 & 9.00 & 9.00 & 22.00 & 8.00 & 8.00 \\
$2^7$ & $2^8$ & 128.00 & 8.00 & 9.00 & 128.00 & 9.00 & 9.00 & 39.00 & 9.00 & 9.00 \\
$2^8$ & $2^8$ & 256.00 & 8.00 & 9.00 & 256.00 & 9.01 & 10.00 & 72.01 & 9.00 & 10.00 \\
\hline
\end{tabular}
\caption{Averaged number of GMRES iterations in forward block solution without preconditioner, with $P_{\alpha,N}$ circulant preconditioner, and $P_{\omega,\alpha,N}$  $\omega-$circulant preconditioner - case $a(x) = x^2+1$.}
\label{tab:iterations_es2_alpha11}
\end{table}
\begin{table}[H]
\centering
\vspace{0.5em}
\small
\begin{tabular}{|cc|ccc|ccc|ccc|}
\hline
& & \multicolumn{9}{c|}{$\alpha=1.5$} \\
& & \multicolumn{3}{c|}{$\beta=0.1$} & \multicolumn{3}{c|}{$\beta=0.5$} & \multicolumn{3}{c|}{$\beta=0.9$} \\
\cline{3-11}
$N$ & $M$
& -  & ${P_{\alpha,N}}$ & $P_{\omega,\alpha,N}$
& -  & ${P_{\alpha,N}}$ & $P_{\omega,\alpha,N}$
& -  & ${P_{\alpha,N}}$ & $P_{\omega,\alpha,N}$
\\
\hline
$2^4$ & $2^4$ & 16.00 & 7.00 & 8.00 & 16.00 & 8.00 & 7.00 & 16.00 & 8.00 & 8.00 \\
$2^5$ & $2^4$ & 32.00 & 7.00 & 8.00 & 32.00 & 8.00 & 7.00 & 32.00 & 8.06 & 8.00 \\
$2^6$ & $2^4$ & 64.00 & 7.00 & 8.00 & 64.00 & 8.00 & 7.00 & 64.00 & 9.00 & 8.94 \\
$2^7$ & $2^4$ & 128.00 & 7.00 & 7.94 & 128.00 & 8.00 & 6.94 & 128.00 & 9.00 & 8.62 \\
$2^8$ & $2^4$ & 256.00 & 7.00 & 7.88 & 256.00 & 8.00 & 6.06 & 256.00 & 8.38 & 8.06 \\
\hline
$2^4$ & $2^5$ & 16.00 & 7.00 & 8.00 & 16.00 & 8.00 & 7.03 & 16.00 & 7.97 & 8.00 \\
$2^5$ & $2^5$ & 32.00 & 6.97 & 8.00 & 32.00 & 8.00 & 8.00 & 32.00 & 8.16 & 9.00 \\
$2^6$ & $2^5$ & 64.00 & 7.00 & 8.00 & 64.00 & 8.03 & 8.00 & 64.00 & 8.97 & 9.00 \\
$2^7$ & $2^5$ & 128.00 & 7.00 & 7.97 & 128.00 & 7.97 & 8.00 & 128.00 & 8.97 & 9.00 \\
$2^8$ & $2^5$ & 256.00 & 7.00 & 7.94 & 256.00 & 8.00 & 7.97 & 256.00 & 8.94 & 9.00 \\
\hline
$2^4$ & $2^6$ & 16.00 & 6.98 & 8.00 & 16.00 & 8.00 & 8.00 & 16.00 & 7.97 & 8.00 \\
$2^5$ & $2^6$ & 32.00 & 6.98 & 8.00 & 32.00 & 8.02 & 8.00 & 32.00 & 8.91 & 9.00 \\
$2^6$ & $2^6$ & 64.00 & 7.00 & 8.00 & 64.00 & 8.11 & 8.00 & 64.00 & 9.00 & 9.00 \\
$2^7$ & $2^6$ & 128.00 & 6.95 & 7.98 & 128.00 & 8.06 & 8.00 & 128.00 & 9.03 & 9.00 \\
$2^8$ & $2^6$ & 256.00 & 7.00 & 7.97 & 256.00 & 7.98 & 8.00 & 256.00 & 9.02 & 9.00 \\
\hline
$2^4$ & $2^7$ & 16.00 & 7.00 & 8.00 & 16.00 & 8.00 & 8.00 & 16.00 & 7.01 & 7.99 \\
$2^5$ & $2^7$ & 32.00 & 7.00 & 8.00 & 32.00 & 8.01 & 8.00 & 32.00 & 8.91 & 9.00 \\
$2^6$ & $2^7$ & 64.00 & 6.99 & 8.00 & 64.00 & 8.99 & 8.99 & 63.00 & 9.00 & 9.00 \\
$2^7$ & $2^7$ & 128.00 & 6.98 & 7.99 & 128.00 & 8.98 & 8.96 & 120.06 & 9.02 & 9.00 \\
$2^8$ & $2^7$ & 256.00 & 6.91 & 7.98 & 256.00 & 8.77 & 8.01 & 231.01 & 9.02 & 9.01 \\
\hline
$2^4$ & $2^8$ & 16.00 & 7.00 & 8.00 & 16.00 & 8.00 & 8.00 & 14.00 & 7.00 & 7.00 \\
$2^5$ & $2^8$ & 32.00 & 7.00 & 8.00 & 32.00 & 8.02 & 9.00 & 25.00 & 8.00 & 9.00 \\
$2^6$ & $2^8$ & 64.00 & 7.01 & 8.00 & 64.00 & 9.00 & 9.00 & 47.00 & 9.00 & 9.00 \\
$2^7$ & $2^8$ & 128.00 & 7.02 & 8.00 & 128.00 & 9.00 & 9.00 & 89.00 & 9.01 & 10.00 \\
$2^8$ & $2^8$ & 256.00 & 6.96 & 7.99 & 256.00 & 8.98 & 9.00 & 170.05 & 9.02 & 10.00 \\
\hline
\end{tabular}
\caption{Averaged number of GMRES iterations in forward block solution without preconditioner, with $P_{\alpha,N}$ circulant preconditioner, and $P_{\omega,\alpha,N}$  $\omega-$circulant preconditioner - case $a(x) = x^2+1$.}
\label{tab:iterations_es2_alpha15}
\end{table}
\begin{table}[H]
\centering
\vspace{0.5em}
\small
\begin{tabular}{|cc|ccc|ccc|ccc|}
\hline
& & \multicolumn{9}{c|}{$\alpha=1.9$} \\
& & \multicolumn{3}{c|}{$\beta=0.1$} & \multicolumn{3}{c|}{$\beta=0.5$} & \multicolumn{3}{c|}{$\beta=0.9$} \\
\cline{3-11}
$N$ & $M$
& -  & ${P_{\alpha,N}}$ & $P_{\omega,\alpha,N}$
& -  & ${P_{\alpha,N}}$ & $P_{\omega,\alpha,N}$
& -  & ${P_{\alpha,N}}$ & $P_{\omega,\alpha,N}$
\\
\hline
$2^4$ & $2^4$ & 16.00 & 6.00 & 7.00 & 16.00 & 7.00 & 7.00 & 16.00 & 7.94 & 7.00 \\
$2^5$ & $2^4$ & 32.00 & 6.00 & 7.00 & 32.00 & 7.00 & 7.00 & 32.00 & 8.00 & 8.00 \\
$2^6$ & $2^4$ & 64.00 & 6.00 & 7.00 & 64.00 & 7.00 & 7.00 & 64.00 & 8.00 & 7.06 \\
$2^7$ & $2^4$ & 128.00 & 6.00 & 6.06 & 128.00 & 6.94 & 6.06 & 128.00 & 7.12 & 7.00 \\
$2^8$ & $2^4$ & 256.00 & 5.00 & 6.00 & 256.00 & 6.00 & 6.00 & 256.00 & 7.00 & 7.00 \\
\hline
$2^4$ & $2^5$ & 16.00 & 6.00 & 7.00 & 16.00 & 7.00 & 6.97 & 16.00 & 8.00 & 8.00 \\
$2^5$ & $2^5$ & 32.00 & 6.00 & 7.00 & 32.00 & 7.00 & 6.00 & 32.00 & 8.03 & 8.00 \\
$2^6$ & $2^5$ & 64.00 & 6.00 & 7.00 & 64.00 & 7.00 & 6.00 & 64.00 & 8.00 & 8.00 \\
$2^7$ & $2^5$ & 128.00 & 6.00 & 6.03 & 128.00 & 7.00 & 5.00 & 128.00 & 8.00 & 8.00 \\
$2^8$ & $2^5$ & 256.00 & 5.00 & 6.00 & 256.00 & 7.00 & 5.00 & 256.00 & 7.97 & 7.97 \\
\hline
$2^4$ & $2^6$ & 16.00 & 6.00 & 7.00 & 16.00 & 7.89 & 7.00 & 16.00 & 8.00 & 8.00 \\
$2^5$ & $2^6$ & 32.00 & 6.00 & 7.00 & 32.00 & 7.95 & 7.00 & 32.00 & 8.91 & 8.00 \\
$2^6$ & $2^6$ & 64.00 & 6.00 & 7.00 & 64.00 & 7.02 & 7.00 & 64.00 & 8.84 & 8.00 \\
$2^7$ & $2^6$ & 128.00 & 6.00 & 6.02 & 128.00 & 7.00 & 7.00 & 128.00 & 8.02 & 8.00 \\
$2^8$ & $2^6$ & 256.00 & 5.77 & 6.00 & 256.00 & 7.00 & 6.00 & 256.00 & 8.00 & 8.00 \\
\hline
$2^4$ & $2^7$ & 16.00 & 6.00 & 7.00 & 16.00 & 7.99 & 7.00 & 16.00 & 8.00 & 8.00 \\
$2^5$ & $2^7$ & 32.00 & 6.00 & 7.00 & 32.00 & 8.00 & 8.00 & 32.00 & 8.99 & 8.99 \\
$2^6$ & $2^7$ & 64.00 & 6.00 & 7.00 & 64.00 & 8.00 & 7.01 & 64.00 & 9.00 & 9.00 \\
$2^7$ & $2^7$ & 128.00 & 6.00 & 6.01 & 128.00 & 7.02 & 7.00 & 128.00 & 8.98 & 8.99 \\
$2^8$ & $2^7$ & 256.00 & 5.97 & 6.01 & 256.00 & 7.00 & 7.00 & 256.00 & 8.01 & 8.00 \\
\hline
$2^4$ & $2^8$ & 16.00 & 6.00 & 7.00 & 16.00 & 8.00 & 8.00 & 16.00 & 7.78 & 8.00 \\
$2^5$ & $2^8$ & 32.00 & 6.00 & 7.00 & 32.00 & 8.00 & 8.00 & 32.00 & 8.64 & 9.00 \\
$2^6$ & $2^8$ & 64.00 & 6.00 & 7.00 & 64.00 & 8.00 & 8.00 & 64.00 & 8.84 & 9.00 \\
$2^7$ & $2^8$ & 128.00 & 6.00 & 6.00 & 128.00 & 8.00 & 8.00 & 128.00 & 9.00 & 9.00 \\
$2^8$ & $2^8$ & 256.00 & 6.00 & 6.00 & 256.00 & 7.00 & 7.00 & 252.08 & 9.00 & 9.00 \\
\hline
\end{tabular}
\caption{Averaged number of GMRES iterations in forward block solution without preconditioner, with $P_{\alpha,N}$ circulant preconditioner, and $P_{\omega,\alpha,N}$  $\omega-$circulant preconditioner - case $a(x) = x^2+1$.}
\label{tab:iterations_es2_alpha19}
\end{table}
%
\begin{table}[H]
\centering
\vspace{0.5em}
\small
\begin{tabular}{|cc|ccc|ccc|ccc|}
\hline
& & \multicolumn{9}{c|}{$\alpha=1.1$} \\
& & \multicolumn{3}{c|}{$\beta=0.1$} & \multicolumn{3}{c|}{$\beta=0.5$} & \multicolumn{3}{c|}{$\beta=0.9$} \\
\cline{3-11}
$N$ & $M$
& -  & ${P_{\alpha,\beta,NM}}$ & $P_{\omega,\alpha,\beta,MN}$
& -  & ${P_{\alpha,\beta,NM}}$ & $P_{\omega,\alpha,\beta,MN}$
& -  & ${P_{\alpha,\beta,NM}}$ & $P_{\omega,\alpha,\beta,MN}$
\\
\hline
$2^4$ & $2^4$ & 33.00 & 11.00 & 10.00 & 40.00 & 11.00 & 10.00 & 46.00 & 10.00 & 9.00 \\
$2^5$ & $2^4$ & 63.00 & 11.00 & 10.00 & 71.00 & 11.00 & 11.00 & 77.00 & 11.00 & 10.00 \\
$2^6$ & $2^4$ & 123.00 & 12.00 & 11.00 & 134.00 & 12.00 & 11.00 & 139.00 & 11.00 & 10.00 \\
$2^7$ & $2^4$ & 240.00 & 13.00 & 11.00 & 261.00 & 13.00 & 11.00 & 263.00 & 12.00 & 11.00 \\
$2^8$ & $2^4$ & 462.00 & 13.00 & 11.00 & 500.00 & 13.00 & 12.00 & 500.00 & 12.00 & 12.00 \\
\hline
$2^4$ & $2^5$ & 33.00 & 11.00 & 10.00 & 42.00 & 11.00 & 10.00 & 60.00 & 10.00 & 10.00 \\
$2^5$ & $2^5$ & 63.00 & 11.00 & 10.00 & 73.00 & 11.00 & 11.00 & 90.00 & 11.00 & 10.00 \\
$2^6$ & $2^5$ & 123.00 & 12.00 & 11.00 & 135.00 & 12.00 & 11.00 & 150.00 & 12.00 & 11.00 \\
$2^7$ & $2^5$ & 240.00 & 13.00 & 11.00 & 262.00 & 13.00 & 11.00 & 271.00 & 12.00 & 11.00 \\
$2^8$ & $2^5$ & 463.00 & 13.00 & 11.00 & 500.00 & 13.00 & 12.00 & 500.00 & 13.00 & 12.00 \\
\hline
$2^4$ & $2^6$ & 33.00 & 11.00 & 10.00 & 44.00 & 11.00 & 10.00 & 86.00 & 10.00 & 10.00 \\
$2^5$ & $2^6$ & 63.00 & 12.00 & 11.00 & 75.00 & 11.00 & 11.00 & 114.00 & 12.00 & 10.00 \\
$2^6$ & $2^6$ & 123.00 & 12.00 & 11.00 & 137.00 & 12.00 & 11.00 & 172.00 & 12.00 & 11.00 \\
$2^7$ & $2^6$ & 240.00 & 13.00 & 11.00 & 263.00 & 13.00 & 12.00 & 291.00 & 13.00 & 12.00 \\
$2^8$ & $2^6$ & 463.00 & 14.00 & 12.00 & 500.00 & 13.00 & 12.00 & 500.00 & 13.00 & 12.00 \\
\hline
$2^4$ & $2^7$ & 33.00 & 12.00 & 10.00 & 47.00 & 11.00 & 10.00 & 140.00 & 11.00 & 10.00 \\
$2^5$ & $2^7$ & 63.00 & 12.00 & 11.00 & 78.00 & 12.00 & 11.00 & 163.00 & 12.00 & 11.00 \\
$2^6$ & $2^7$ & 123.00 & 12.00 & 11.00 & 139.00 & 12.00 & 11.00 & 218.00 & 12.00 & 11.00 \\
$2^7$ & $2^7$ & 240.00 & 13.00 & 11.00 & 265.00 & 13.00 & 12.00 & 332.00 & 13.00 & 12.00 \\
$2^8$ & $2^7$ & 463.00 & 14.00 & 12.00 & 500.00 & 14.00 & 12.00 & 500.00 & 13.00 & 12.00 \\
\hline
$2^4$ & $2^8$ & 33.00 & 12.00 & 10.00 & 51.00 & 11.00 & 10.00 & 238.00 & 11.00 & 10.00 \\
$2^5$ & $2^8$ & 63.00 & 12.00 & 11.00 & 81.00 & 12.00 & 11.00 & 259.00 & 12.00 & 11.00 \\
$2^6$ & $2^8$ & 123.00 & 12.00 & 11.00 & 143.00 & 12.00 & 11.00 & 310.00 & 12.00 & 11.00 \\
$2^7$ & $2^8$ & 240.00 & 13.00 & 12.00 & 268.00 & 13.00 & 12.00 & 420.00 & 13.00 & 12.00 \\
$2^8$ & $2^8$ & 463.00 & 14.00 & 12.00 & 500.00 & 14.00 & 12.00 & 500.00 & 13.00 & 12.00 \\
\hline
\end{tabular}
\caption{Number of GMRES iterations without preconditioner, with $P_{\alpha,\beta,NM}$ circulant preconditioner, and $P_{\omega,\alpha,\beta,MN}$  $\omega-$circulant preconditioner - case $a(x) = 1$.}
\label{tab:iterations_tot_es1_alpha11}
\end{table}
\begin{table}[H]
\centering
\vspace{0.5em}
\small
\begin{tabular}{|cc|ccc|ccc|ccc|}
\hline
& & \multicolumn{9}{c|}{$\alpha=1.5$} \\
& & \multicolumn{3}{c|}{$\beta=0.1$} & \multicolumn{3}{c|}{$\beta=0.5$} & \multicolumn{3}{c|}{$\beta=0.9$} \\
\cline{3-11}
$N$ & $M$
& -  & ${P_{\alpha,\beta,NM}}$ & $P_{\omega,\alpha,\beta,MN}$
& -  & ${P_{\alpha,\beta,NM}}$ & $P_{\omega,\alpha,\beta,MN}$
& -  & ${P_{\alpha,\beta,NM}}$ & $P_{\omega,\alpha,\beta,MN}$
\\
\hline
$2^4$ & $2^4$ & 31.00 & 13.00 & 13.00 & 46.00 & 17.00 & 16.00 & 54.00 & 15.00 & 15.00 \\
$2^5$ & $2^4$ & 58.00 & 14.00 & 13.00 & 79.00 & 18.00 & 18.00 & 91.00 & 17.00 & 17.00 \\
$2^6$ & $2^4$ & 111.00 & 14.00 & 14.00 & 144.00 & 19.00 & 19.00 & 164.00 & 18.00 & 18.00 \\
$2^7$ & $2^4$ & 218.00 & 15.00 & 14.00 & 272.00 & 20.00 & 18.00 & 307.00 & 19.00 & 19.00 \\
$2^8$ & $2^4$ & 416.00 & 15.00 & 15.00 & 500.00 & 20.00 & 18.00 & 500.00 & 20.00 & 20.00 \\
\hline
$2^4$ & $2^5$ & 31.00 & 14.00 & 13.00 & 48.00 & 19.00 & 17.00 & 70.00 & 16.00 & 15.00 \\
$2^5$ & $2^5$ & 58.00 & 15.00 & 14.00 & 82.00 & 20.00 & 18.00 & 114.00 & 18.00 & 18.00 \\
$2^6$ & $2^5$ & 111.00 & 15.00 & 14.00 & 147.00 & 20.00 & 19.00 & 204.00 & 19.00 & 20.00 \\
$2^7$ & $2^5$ & 218.00 & 16.00 & 15.00 & 276.00 & 22.00 & 19.00 & 379.00 & 20.00 & 21.00 \\
$2^8$ & $2^5$ & 417.00 & 16.00 & 16.00 & 500.00 & 23.00 & 20.00 & 500.00 & 22.00 & 21.00 \\
\hline
$2^4$ & $2^6$ & 31.00 & 14.00 & 13.00 & 50.00 & 19.00 & 17.00 & 100.00 & 16.00 & 16.00 \\
$2^5$ & $2^6$ & 58.00 & 15.00 & 14.00 & 84.00 & 20.00 & 19.00 & 148.00 & 19.00 & 18.00 \\
$2^6$ & $2^6$ & 111.00 & 16.00 & 15.00 & 150.00 & 21.00 & 20.00 & 250.00 & 20.00 & 20.00 \\
$2^7$ & $2^6$ & 218.00 & 16.00 & 15.00 & 278.00 & 23.00 & 21.00 & 457.00 & 22.00 & 21.00 \\
$2^8$ & $2^6$ & 418.00 & 17.00 & 16.00 & 500.00 & 25.00 & 22.00 & 500.00 & 24.00 & 23.00 \\
\hline
$2^4$ & $2^7$ & 31.00 & 15.00 & 14.00 & 52.00 & 19.00 & 17.00 & 153.00 & 17.00 & 16.00 \\
$2^5$ & $2^7$ & 58.00 & 15.00 & 14.00 & 86.00 & 21.00 & 19.00 & 208.00 & 19.00 & 19.00 \\
$2^6$ & $2^7$ & 111.00 & 16.00 & 15.00 & 153.00 & 22.00 & 21.00 & 329.00 & 21.00 & 21.00 \\
$2^7$ & $2^7$ & 218.00 & 17.00 & 16.00 & 280.00 & 24.00 & 22.00 & 500.00 & 24.00 & 23.00 \\
$2^8$ & $2^7$ & 418.00 & 18.00 & 17.00 & 500.00 & 26.00 & 23.00 & 500.00 & 25.00 & 25.00 \\
\hline
$2^4$ & $2^8$ & 31.00 & 15.00 & 14.00 & 55.00 & 20.00 & 18.00 & 244.00 & 17.00 & 16.00 \\
$2^5$ & $2^8$ & 58.00 & 16.00 & 15.00 & 88.00 & 21.00 & 20.00 & 303.00 & 19.00 & 19.00 \\
$2^6$ & $2^8$ & 111.00 & 17.00 & 16.00 & 155.00 & 23.00 & 21.00 & 443.00 & 22.00 & 21.00 \\
$2^7$ & $2^8$ & 218.00 & 17.00 & 16.00 & 284.00 & 25.00 & 23.00 & 500.00 & 25.00 & 24.00 \\
$2^8$ & $2^8$ & 418.00 & 18.00 & 17.00 & 500.00 & 28.00 & 24.00 & 500.00 & 27.00 & 26.00 \\
\hline
\end{tabular}
\caption{Number of GMRES iterations without preconditioner, with $P_{\alpha,\beta,NM}$ circulant preconditioner, and $P_{\omega,\alpha,\beta,MN}$  $\omega-$circulant preconditioner - case $a(x) = 1$.}
\label{tab:iterations_tot_es1_alpha15}
\end{table}
\begin{table}[H]
\centering
\vspace{0.5em}
\small
\begin{tabular}{|cc|ccc|ccc|ccc|}
\hline
& & \multicolumn{9}{c|}{$\alpha=1.9$} \\
& & \multicolumn{3}{c|}{$\beta=0.1$} & \multicolumn{3}{c|}{$\beta=0.5$} & \multicolumn{3}{c|}{$\beta=0.9$} \\
\cline{3-11}
$N$ & $M$
& -  & ${P_{\alpha,\beta,NM}}$ & $P_{\omega,\alpha,\beta,MN}$
& -  & ${P_{\alpha,\beta,NM}}$ & $P_{\omega,\alpha,\beta,MN}$
& -  & ${P_{\alpha,\beta,NM}}$ & $P_{\omega,\alpha,\beta,MN}$
\\
\hline
$2^4$ & $2^4$ & 36.00 & 16.00 & 13.00 & 56.00 & 21.00 & 19.00 & 78.00 & 18.00 & 18.00 \\
$2^5$ & $2^4$ & 66.00 & 17.00 & 14.00 & 99.00 & 20.00 & 20.00 & 136.00 & 19.00 & 21.00 \\
$2^6$ & $2^4$ & 125.00 & 18.00 & 15.00 & 175.00 & 20.00 & 20.00 & 243.00 & 20.00 & 20.00 \\
$2^7$ & $2^4$ & 219.00 & 19.00 & 15.00 & 312.00 & 20.00 & 20.00 & 437.00 & 20.00 & 20.00 \\
$2^8$ & $2^4$ & 410.00 & 19.00 & 16.00 & 500.00 & 20.00 & 20.00 & 500.00 & 20.00 & 20.00 \\
\hline
$2^4$ & $2^5$ & 35.00 & 17.00 & 14.00 & 59.00 & 23.00 & 22.00 & 99.00 & 19.00 & 20.00 \\
$2^5$ & $2^5$ & 66.00 & 18.00 & 15.00 & 101.00 & 25.00 & 23.00 & 171.00 & 22.00 & 22.00 \\
$2^6$ & $2^5$ & 125.00 & 19.00 & 16.00 & 186.00 & 26.00 & 25.00 & 299.00 & 24.00 & 24.00 \\
$2^7$ & $2^5$ & 219.00 & 20.00 & 17.00 & 316.00 & 27.00 & 27.00 & 500.00 & 26.00 & 26.00 \\
$2^8$ & $2^5$ & 408.00 & 21.00 & 17.00 & 500.00 & 29.00 & 28.00 & 500.00 & 28.00 & 28.00 \\
\hline
$2^4$ & $2^6$ & 35.00 & 18.00 & 15.00 & 61.00 & 25.00 & 22.00 & 130.00 & 20.00 & 20.00 \\
$2^5$ & $2^6$ & 66.00 & 19.00 & 16.00 & 105.00 & 27.00 & 25.00 & 215.00 & 23.00 & 23.00 \\
$2^6$ & $2^6$ & 125.00 & 20.00 & 17.00 & 187.00 & 29.00 & 28.00 & 374.00 & 27.00 & 27.00 \\
$2^7$ & $2^6$ & 219.00 & 21.00 & 17.00 & 322.00 & 31.00 & 30.00 & 500.00 & 30.00 & 30.00 \\
$2^8$ & $2^6$ & 408.00 & 23.00 & 18.00 & 500.00 & 33.00 & 32.00 & 500.00 & 32.00 & 32.00 \\
\hline
$2^4$ & $2^7$ & 35.00 & 19.00 & 15.00 & 63.00 & 26.00 & 23.00 & 178.00 & 21.00 & 21.00 \\
$2^5$ & $2^7$ & 66.00 & 20.00 & 16.00 & 108.00 & 29.00 & 26.00 & 280.00 & 25.00 & 24.00 \\
$2^6$ & $2^7$ & 125.00 & 21.00 & 17.00 & 189.00 & 32.00 & 30.00 & 476.00 & 30.00 & 29.00 \\
$2^7$ & $2^7$ & 219.00 & 22.00 & 18.00 & 326.00 & 34.00 & 33.00 & 500.00 & 33.00 & 34.00 \\
$2^8$ & $2^7$ & 407.00 & 24.00 & 19.00 & 500.00 & 37.00 & 35.00 & 500.00 & 36.00 & 37.00 \\
\hline
$2^4$ & $2^8$ & 35.00 & 19.00 & 16.00 & 64.00 & 26.00 & 24.00 & 253.00 & 21.00 & 21.00 \\
$2^5$ & $2^8$ & 66.00 & 21.00 & 17.00 & 110.00 & 31.00 & 27.00 & 374.00 & 26.00 & 26.00 \\
$2^6$ & $2^8$ & 125.00 & 22.00 & 18.00 & 192.00 & 34.00 & 32.00 & 500.00 & 30.00 & 30.00 \\
$2^7$ & $2^8$ & 219.00 & 23.00 & 19.00 & 329.00 & 37.00 & 35.00 & 500.00 & 33.00 & 36.00 \\
$2^8$ & $2^8$ & 407.00 & 24.00 & 20.00 & 500.00 & 39.00 & 37.00 & 500.00 & 36.00 & 41.00 \\
\hline
\end{tabular}
\caption{Number of GMRES iterations without preconditioner, with $P_{\alpha,\beta,NM}$ circulant preconditioner, and $P_{\omega,\alpha,\beta,MN}$  $\omega-$circulant preconditioner - case $a(x) = 1$.}
\label{tab:iterations_tot_es1_alpha19}
\end{table}
\begin{table}[H]
\centering
\vspace{0.5em}
\small
\begin{tabular}{|cc|ccc|ccc|ccc|}
\hline
& & \multicolumn{9}{c|}{$\alpha=1.1$} \\
& & \multicolumn{3}{c|}{$\beta=0.1$} & \multicolumn{3}{c|}{$\beta=0.5$} & \multicolumn{3}{c|}{$\beta=0.9$} \\
\cline{3-11}
$N$ & $M$
& -  & ${P_{\alpha,\beta,NM}}$ & $P_{\omega,\alpha,\beta,MN}$
& -  & ${P_{\alpha,\beta,NM}}$ & $P_{\omega,\alpha,\beta,MN}$
& -  & ${P_{\alpha,\beta,NM}}$ & $P_{\omega,\alpha,\beta,MN}$
\\
\hline
$2^4$ & $2^4$ & 33.00 & 13.00 & 13.00 & 40.00 & 13.00 & 12.00 & 48.00 & 12.00 & 11.00 \\
$2^5$ & $2^4$ & 63.00 & 14.00 & 14.00 & 71.00 & 13.00 & 13.00 & 79.00 & 12.00 & 12.00 \\
$2^6$ & $2^4$ & 123.00 & 15.00 & 14.00 & 134.00 & 14.00 & 13.00 & 142.00 & 13.00 & 12.00 \\
$2^7$ & $2^4$ & 238.00 & 15.00 & 14.00 & 264.00 & 14.00 & 13.00 & 270.00 & 13.00 & 13.00 \\
$2^8$ & $2^4$ & 475.00 & 15.00 & 14.00 & 500.00 & 15.00 & 14.00 & 500.00 & 13.00 & 13.00 \\
\hline
$2^4$ & $2^5$ & 33.00 & 14.00 & 13.00 & 41.00 & 13.00 & 12.00 & 62.00 & 12.00 & 11.00 \\
$2^5$ & $2^5$ & 63.00 & 14.00 & 14.00 & 73.00 & 14.00 & 13.00 & 93.00 & 13.00 & 12.00 \\
$2^6$ & $2^5$ & 123.00 & 15.00 & 14.00 & 135.00 & 14.00 & 13.00 & 156.00 & 13.00 & 13.00 \\
$2^7$ & $2^5$ & 239.00 & 15.00 & 15.00 & 265.00 & 15.00 & 14.00 & 284.00 & 14.00 & 13.00 \\
$2^8$ & $2^5$ & 476.00 & 16.00 & 15.00 & 500.00 & 15.00 & 14.00 & 500.00 & 14.00 & 13.00 \\
\hline
$2^4$ & $2^6$ & 33.00 & 14.00 & 13.00 & 43.00 & 13.00 & 12.00 & 88.00 & 13.00 & 12.00 \\
$2^5$ & $2^6$ & 63.00 & 14.00 & 14.00 & 74.00 & 14.00 & 13.00 & 117.00 & 14.00 & 12.00 \\
$2^6$ & $2^6$ & 123.00 & 15.00 & 15.00 & 137.00 & 14.00 & 13.00 & 179.00 & 14.00 & 13.00 \\
$2^7$ & $2^6$ & 239.00 & 16.00 & 15.00 & 266.00 & 15.00 & 14.00 & 306.00 & 14.00 & 13.00 \\
$2^8$ & $2^6$ & 476.00 & 16.00 & 15.00 & 500.00 & 15.00 & 14.00 & 500.00 & 15.00 & 14.00 \\
\hline
$2^4$ & $2^7$ & 33.00 & 14.00 & 13.00 & 46.00 & 14.00 & 13.00 & 138.00 & 13.00 & 12.00 \\
$2^5$ & $2^7$ & 63.00 & 15.00 & 14.00 & 77.00 & 14.00 & 13.00 & 165.00 & 14.00 & 13.00 \\
$2^6$ & $2^7$ & 123.00 & 15.00 & 15.00 & 139.00 & 14.00 & 13.00 & 225.00 & 14.00 & 13.00 \\
$2^7$ & $2^7$ & 239.00 & 16.00 & 15.00 & 268.00 & 15.00 & 14.00 & 351.00 & 15.00 & 14.00 \\
$2^8$ & $2^7$ & 476.00 & 16.00 & 15.00 & 500.00 & 16.00 & 14.00 & 500.00 & 15.00 & 14.00 \\
\hline
$2^4$ & $2^8$ & 33.00 & 14.00 & 13.00 & 49.00 & 14.00 & 13.00 & 231.00 & 13.00 & 12.00 \\
$2^5$ & $2^8$ & 63.00 & 15.00 & 14.00 & 79.00 & 14.00 & 13.00 & 255.00 & 14.00 & 13.00 \\
$2^6$ & $2^8$ & 123.00 & 15.00 & 15.00 & 141.00 & 14.00 & 13.00 & 312.00 & 15.00 & 14.00 \\
$2^7$ & $2^8$ & 239.00 & 16.00 & 15.00 & 269.00 & 15.00 & 14.00 & 435.00 & 15.00 & 14.00 \\
$2^8$ & $2^8$ & 476.00 & 17.00 & 16.00 & 500.00 & 16.00 & 14.00 & 500.00 & 16.00 & 15.00 \\
\hline
\end{tabular}
\caption{Number of GMRES iterations without preconditioner, with $P_{\alpha,\beta,NM}$ circulant preconditioner, and $P_{\omega,\alpha,\beta,MN}$  $\omega-$circulant preconditioner - case $a(x) = x^2+1$.}
\label{tab:iterations_tot_es2_alpha11}
\end{table}
\begin{table}[H]
\centering
\vspace{0.5em}
\small
\begin{tabular}{|cc|ccc|ccc|ccc|}
\hline
& & \multicolumn{9}{c|}{$\alpha=1.5$} \\
& & \multicolumn{3}{c|}{$\beta=0.1$} & \multicolumn{3}{c|}{$\beta=0.5$} & \multicolumn{3}{c|}{$\beta=0.9$} \\
\cline{3-11}
$N$ & $M$
& -  & ${P_{\alpha,\beta,NM}}$ & $P_{\omega,\alpha,\beta,MN}$
& -  & ${P_{\alpha,\beta,NM}}$ & $P_{\omega,\alpha,\beta,MN}$
& -  & ${P_{\alpha,\beta,NM}}$ & $P_{\omega,\alpha,\beta,MN}$
\\
\hline
$2^4$ & $2^4$ & 33.00 & 16.00 & 16.00 & 49.00 & 19.00 & 18.00 & 61.00 & 17.00 & 16.00 \\
$2^5$ & $2^4$ & 62.00 & 18.00 & 17.00 & 86.00 & 20.00 & 19.00 & 103.00 & 18.00 & 18.00 \\
$2^6$ & $2^4$ & 119.00 & 18.00 & 18.00 & 157.00 & 21.00 & 20.00 & 186.00 & 19.00 & 20.00 \\
$2^7$ & $2^4$ & 230.00 & 19.00 & 18.00 & 290.00 & 22.00 & 21.00 & 345.00 & 20.00 & 21.00 \\
$2^8$ & $2^4$ & 402.00 & 20.00 & 19.00 & 500.00 & 22.00 & 21.00 & 500.00 & 22.00 & 22.00 \\
\hline
$2^4$ & $2^5$ & 33.00 & 16.00 & 16.00 & 52.00 & 20.00 & 19.00 & 78.00 & 17.00 & 17.00 \\
$2^5$ & $2^5$ & 62.00 & 18.00 & 18.00 & 88.00 & 22.00 & 21.00 & 129.00 & 19.00 & 19.00 \\
$2^6$ & $2^5$ & 119.00 & 19.00 & 19.00 & 160.00 & 23.00 & 22.00 & 231.00 & 21.00 & 21.00 \\
$2^7$ & $2^5$ & 231.00 & 20.00 & 19.00 & 298.00 & 24.00 & 23.00 & 430.00 & 23.00 & 22.00 \\
$2^8$ & $2^5$ & 403.00 & 21.00 & 20.00 & 500.00 & 26.00 & 24.00 & 500.00 & 24.00 & 23.00 \\
\hline
$2^4$ & $2^6$ & 33.00 & 17.00 & 16.00 & 53.00 & 21.00 & 19.00 & 108.00 & 18.00 & 17.00 \\
$2^5$ & $2^6$ & 62.00 & 19.00 & 18.00 & 90.00 & 23.00 & 22.00 & 166.00 & 20.00 & 19.00 \\
$2^6$ & $2^6$ & 119.00 & 20.00 & 19.00 & 162.00 & 25.00 & 23.00 & 284.00 & 23.00 & 22.00 \\
$2^7$ & $2^6$ & 231.00 & 21.00 & 20.00 & 302.00 & 26.00 & 25.00 & 500.00 & 25.00 & 24.00 \\
$2^8$ & $2^6$ & 403.00 & 22.00 & 21.00 & 500.00 & 28.00 & 26.00 & 500.00 & 27.00 & 26.00 \\
\hline
$2^4$ & $2^7$ & 33.00 & 17.00 & 16.00 & 55.00 & 21.00 & 20.00 & 159.00 & 18.00 & 17.00 \\
$2^5$ & $2^7$ & 62.00 & 19.00 & 19.00 & 92.00 & 24.00 & 22.00 & 226.00 & 21.00 & 20.00 \\
$2^6$ & $2^7$ & 119.00 & 21.00 & 20.00 & 164.00 & 26.00 & 24.00 & 365.00 & 23.00 & 22.00 \\
$2^7$ & $2^7$ & 231.00 & 22.00 & 21.00 & 304.00 & 27.00 & 26.00 & 500.00 & 26.00 & 25.00 \\
$2^8$ & $2^7$ & 404.00 & 23.00 & 22.00 & 500.00 & 29.00 & 27.00 & 500.00 & 29.00 & 28.00 \\
\hline
$2^4$ & $2^8$ & 33.00 & 18.00 & 17.00 & 57.00 & 22.00 & 20.00 & 243.00 & 18.00 & 17.00 \\
$2^5$ & $2^8$ & 62.00 & 19.00 & 19.00 & 94.00 & 24.00 & 22.00 & 317.00 & 21.00 & 20.00 \\
$2^6$ & $2^8$ & 119.00 & 21.00 & 20.00 & 166.00 & 27.00 & 25.00 & 481.00 & 24.00 & 23.00 \\
$2^7$ & $2^8$ & 231.00 & 22.00 & 21.00 & 306.00 & 28.00 & 27.00 & 500.00 & 28.00 & 27.00 \\
$2^8$ & $2^8$ & 404.00 & 23.00 & 22.00 & 500.00 & 30.00 & 29.00 & 500.00 & 31.00 & 30.00 \\
\hline
\end{tabular}
\caption{Number of GMRES iterations without preconditioner, with $P_{\alpha,\beta,NM}$ circulant preconditioner, and $P_{\omega,\alpha,\beta,MN}$  $\omega-$circulant preconditioner - case $a(x) = x^2+1$.}
\label{tab:iterations_tot_es2_alpha15}
\end{table}
\begin{table}[H]
\centering
\vspace{0.5em}
\small
\begin{tabular}{|cc|ccc|ccc|ccc|}
\hline
& & \multicolumn{9}{c|}{$\alpha=1.9$} \\
& & \multicolumn{3}{c|}{$\beta=0.1$} & \multicolumn{3}{c|}{$\beta=0.5$} & \multicolumn{3}{c|}{$\beta=0.9$} \\
\cline{3-11}
$N$ & $M$
& -  & ${P_{\alpha,\beta,NM}}$ & $P_{\omega,\alpha,\beta,MN}$
& -  & ${P_{\alpha,\beta,NM}}$ & $P_{\omega,\alpha,\beta,MN}$
& -  & ${P_{\alpha,\beta,NM}}$ & $P_{\omega,\alpha,\beta,MN}$
\\
\hline
$2^4$ & $2^4$ & 38.00 & 19.00 & 18.00 & 59.00 & 21.00 & 21.00 & 85.00 & 19.00 & 20.00 \\
$2^5$ & $2^4$ & 74.00 & 20.00 & 19.00 & 105.00 & 21.00 & 21.00 & 151.00 & 20.00 & 21.00 \\
$2^6$ & $2^4$ & 132.00 & 20.00 & 20.00 & 193.00 & 21.00 & 21.00 & 272.00 & 21.00 & 21.00 \\
$2^7$ & $2^4$ & 253.00 & 19.00 & 19.00 & 340.00 & 21.00 & 21.00 & 488.00 & 21.00 & 21.00 \\
$2^8$ & $2^4$ & 471.00 & 20.00 & 19.00 & 500.00 & 20.00 & 20.00 & 500.00 & 20.00 & 21.00 \\
\hline
$2^4$ & $2^5$ & 38.00 & 20.00 & 19.00 & 61.00 & 25.00 & 25.00 & 106.00 & 21.00 & 22.00 \\
$2^5$ & $2^5$ & 74.00 & 21.00 & 20.00 & 111.00 & 28.00 & 27.00 & 186.00 & 24.00 & 25.00 \\
$2^6$ & $2^5$ & 132.00 & 22.00 & 22.00 & 196.00 & 30.00 & 29.00 & 330.00 & 26.00 & 27.00 \\
$2^7$ & $2^5$ & 252.00 & 21.00 & 20.00 & 347.00 & 31.00 & 30.00 & 500.00 & 28.00 & 29.00 \\
$2^8$ & $2^5$ & 465.00 & 22.00 & 21.00 & 500.00 & 30.00 & 29.00 & 500.00 & 29.00 & 31.00 \\
\hline
$2^4$ & $2^6$ & 38.00 & 21.00 & 20.00 & 64.00 & 27.00 & 26.00 & 135.00 & 23.00 & 23.00 \\
$2^5$ & $2^6$ & 74.00 & 22.00 & 21.00 & 112.00 & 31.00 & 30.00 & 231.00 & 26.00 & 26.00 \\
$2^6$ & $2^6$ & 133.00 & 24.00 & 23.00 & 201.00 & 34.00 & 33.00 & 407.00 & 29.00 & 30.00 \\
$2^7$ & $2^6$ & 252.00 & 23.00 & 22.00 & 348.00 & 36.00 & 34.00 & 500.00 & 32.00 & 34.00 \\
$2^8$ & $2^6$ & 461.00 & 23.00 & 22.00 & 500.00 & 35.00 & 34.00 & 500.00 & 34.00 & 36.00 \\
\hline
$2^4$ & $2^7$ & 38.00 & 22.00 & 21.00 & 65.00 & 29.00 & 28.00 & 178.00 & 23.00 & 24.00 \\
$2^5$ & $2^7$ & 74.00 & 23.00 & 22.00 & 115.00 & 33.00 & 31.00 & 294.00 & 27.00 & 28.00 \\
$2^6$ & $2^7$ & 133.00 & 25.00 & 24.00 & 205.00 & 37.00 & 35.00 & 500.00 & 31.00 & 32.00 \\
$2^7$ & $2^7$ & 252.00 & 24.00 & 23.00 & 350.00 & 39.00 & 36.00 & 500.00 & 35.00 & 37.00 \\
$2^8$ & $2^7$ & 459.00 & 24.00 & 23.00 & 500.00 & 39.00 & 38.00 & 500.00 & 36.00 & 40.00 \\
\hline
$2^4$ & $2^8$ & 38.00 & 23.00 & 21.00 & 67.00 & 31.00 & 29.00 & 244.00 & 23.00 & 24.00 \\
$2^5$ & $2^8$ & 74.00 & 24.00 & 23.00 & 117.00 & 34.00 & 33.00 & 384.00 & 27.00 & 29.00 \\
$2^6$ & $2^8$ & 133.00 & 26.00 & 25.00 & 207.00 & 38.00 & 36.00 & 500.00 & 31.00 & 34.00 \\
$2^7$ & $2^8$ & 252.00 & 24.00 & 24.00 & 358.00 & 40.00 & 39.00 & 500.00 & 35.00 & 39.00 \\
$2^8$ & $2^8$ & 458.00 & 25.00 & 24.00 & 500.00 & 42.00 & 41.00 & 500.00 & 37.00 & 43.00 \\
\hline
\end{tabular}
\caption{Number of GMRES iterations without preconditioner, with $P_{\alpha,\beta,NM}$ circulant preconditioner, and $P_{\omega,\alpha,\beta,MN}$  $\omega-$circulant preconditioner - case $a(x) = x^2+1$.}
\label{tab:iterations_tot_es2_alpha19}
\end{table}

\section{Conclusion}\label{conclusion}
In this paper, we investigated all-at-once linear systems arising from variable-coefficient fractional evolution equations with weakly singular temporal kernels and Riemann--Liouville space fractional derivatives. The temporal discretization was carried out using the L1 scheme, while the spatial fractional operator was approximated through finite difference techniques, leading to large structured $(d+1)$-level Toeplitz-like systems, where one level corresponds to the temporal discretization and the remaining $d$ levels arise from the spatial discretization. In this work, we focused on the case $d = 1$, resulting in a two-level block structured system. Using the GLT framework, we analyzed the asymptotic spectral distribution of the resulting matrix sequences and derived the associated GLT symbols for both the spatial and temporal discretization operators.

Motivated by the structure of the coefficient matrix, we proposed a block lower triangular preconditioner obtained by simplifying the upper triangular part of the spatial operator while preserving the main spectral properties of the original system. The theoretical analysis showed that the preconditioned matrix sequence are well clustered in the sense of the singular values

The numerical experiments confirmed the effectiveness of the proposed preconditioning strategy for both constant and variable diffusion coefficients. In all tested cases, the preconditioners significantly reduced the condition numbers of the discretized systems and improved the convergence behavior of the GMRES method by lowering the number of iterations and computational cost. The observed spectral clustering was fully consistent with the theoretical GLT analysis.

\section*{Acknowledgments}
The research of Muhammad Faisal Khan, Department of Numerical Mathematics, Charles University, Prague, Czech Republic, is supported by the PRIMUS grant PRIMUS/25/SCI/022 of Charles University. The authors are indebted with Cristina Tablino Possio for her contribution in the preconditioning proposals and in the numerical experiments.

\end{document}